\documentclass{amsart}
\newtheorem{defn}{\rm DEFINITION}[section]
\newtheorem{lem}[defn]{\rm LEMMA}
\newtheorem{thm}[defn]{\rm THEOREM}
\newtheorem{prop}[defn]{\rm PROPOSITION}
\newtheorem{cor}[defn]{\rm COROLLARY}
\newtheorem{rem}[defn]{\rm REMARK}
\newtheorem{ex}[defn]{\rm EXAMPLE}

\newcommand{\pp}{{\mathbb P}}
\newcommand{\C}{{\mathbb C}}
\newcommand{\p}{{\mathbb P}}
\newcommand{\pas}{\par\noindent}

\newcommand{\map}{\dasharrow}

\begin{document}
\large
\setcounter{page} {1}
\def\qed{\hfill\vrule height5true pt width5true pt depth0true pt \medskip}
\vfill\eject

\title{ On the classification of defective threefolds}
\date{}

\thispagestyle{empty}

\author{Luca Chiantini}
\address{Dipartimento di Matematica\\
Universit\'a di Siena\\ Via del Capitano, 15\\
 53100 Siena, Italia}
\email{chiantini@unisi.it}

\author{Ciro Ciliberto}
\address{Dipartimento di Matematica\\
Universit\'a di Roma Tor Vergata\\ Via della Ricerca Scientifica\\
 00133 Roma, Italia}
\email{cilibert@axp.mat.uniroma2.it}

\begin{abstract}
In this paper we give the full classification of irreducible projective threefolds
whose $k$-secant variety has dimension smaller than the expected, for some
$k\geq 2$ (see theorem \ref{main} below). As pointed out in the introduction, 
the case $k=1$ was already known before. 
\end{abstract}

\maketitle

\section* {INTRODUCTION} \label {intro}

An irreducible, non--degenerate, projective variety $X$ of dimension
$n$ in $\pp^r$ is called $k$--{\it defective} if its $k$--secant variety $S^k(X)$ has
dimension $s^{(k)}(X)$ smaller than the expected $\sigma^{(k)}(X)$,
which is the minimum between $r$ and $n(k+1)+k$. The difference
$\delta_k(X)=\sigma^{(k)}(X)-s^{(k)}(X)$ is called the $k$--{\it defect}
of $X$.\par

The classification of defective varieties is important in the study of
projective geometry and its applications. The subject goes back to several classical
authors, like Terracini (\cite{Terr1}), Palatini (\cite{Pal}) and
Scorza (\cite{Scorza}), to mention a few. More recently the interest on defective
varieties has been renewed by Zak's spectacular results on the classification of some important classes of smooth varieties (e.g. Severi varieties, Scorza varieties etc., \cite {Zak}).\par

Beside the intrinsic interest of the subject, it turns out that the
understanding of defective varieties is relevant also in other fields of Mathematics:
expressions of polynomials as sums of powers and
Waring type problems, polynomial interpolation,  rank tensor computations and
canonical forms, Bayesian networks, algebraic statistics etc. (see
\cite{Ci2} as a general reference, \cite{CGG}, \cite
{GSS}, \cite {IK}). Se also \cite {So} and \cite{Voisin} as further
examples of results on the subject which have nice applications to number
theoretic problems.\par

In the present paper, we give the complete classification of complex defective threefolds.
\par

One knows that curves are never defective (see \cite{Dale1}). Defective
surfaces have been classified by Terracini in \cite {Terr2}. Terracini's result has been revisited and extended in \cite {Dale2}, \cite {CJ} and \cite
{WDV}. The case of $1$--defective threefolds goes back to Scorza \cite
{Scorza} and has been revisited in  \cite {ChCi}. The case of smooth $1$--defective threefolds was also examined by Fujita \cite {Fuj} and Fujita--Roberts  \cite
{FujRob}. Here we deal with the case of $k$--defective
threefolds, for any $k>1$, and we classify {\it minimally}
$k$--defective threefolds, i.e. $k$--defective threefolds that are not $h$--defective for any $h<k$. This is of course sufficient to describe all defective threefolds. \par

The classification is obtained by using the classical tool of tangential
projections. The basic invariants are the {\it tangential contact
loci} (see \S \ref {notaz}). An important tool is also provided by
Castelnuovo's theory on the growth of Hilbert functions (see \cite {WDV}, \S 6). After having recalled several basic notation, definitions and results in \S \ref {notaz}, we prove some generalities on defective threefolds in \S \ref
{defects} and on the contact loci in \S \ref {2?}. Then we prove our
classification results in \S\S \ref {irrediv}, \ref {red}, \ref {irredcurv}. Indeed \S \ref {irrediv} is devoted to the analysis of the case in which the appropriate tangential contact locus is an irreducible divisor, \S \ref {red} to the case in which the contact locus is reducible, \S \ref {irredcurv} to the
case in which the contact locus is an irreducible curve. \par

Now we state our classification result, obtained by summarizing Theorems
\ref {dimY3}, \ref {dimY2}, \ref {reducib},
\ref {reducib} and various corollaries and remarks accompanying them (for
the definition of the invariant $n_k(X)$ see \ref
{tangentialproj} below): 

\begin{thm} \label {main} Let $X\subset \pp^r$ be an irreducible,
non--degenerate, projective, minimally $k$--defective threefold with $k\geq 2$. Then $X$
is in the following list:

\begin{itemize}

\item [(1)]  $X$ is contained in a cone over the $2$--uple embedding of
a threefold $Y$ of minimal degree $k-1$ in $\pp^{k+1}$, with vertex either a
point, hence $r=4k+2$ and $\delta_k(X)=1$, $s^{(k)}(X)=r-1=4k+1$,
$n_k(X)=1$, or a line, hence $r=4k+3$ and $\delta_k(X)=1$,
$s^{(k)}(X)=r-1=4k+2$, $n_k(X)=2$ (see Example \ref {exx}, (1));\pas

\item [(2)] $k=3$ and either $r=4k+2=14$, $s^{(k)}(X)=r-1=13$,
$n_k(X)=1$, and $X$ is the $2$--uple embedding of a hypersurface $Y$ in $\pp^4$
with $\deg(Y)\geq 3$ or $r=4k+3=15$, $s^{(k)}(X)=r-1=14$,
$n_k(X)=2$, and $X$ is contained in the cone with vertex a point over the
$2$--uple embedding of a hypersurface $Y$ as above (see Example \ref {exx},
(2));\pas

\item [(3)] either $r=4k+2$, $\delta_k(X)=2$, $n_k(X)=1$,
$s^{(k)}(Y')=r-1=4k+1$, and $X$ is the $2$--uple embedding of a threefold $Y$ of degree $k$ in $\pp^{k+1}$ with curve sections of arithmetic
genus $1$ or $r=4k+3$, $\delta_k(X)=1$, $n_k(X)=2$, $s^{(k)}(Y')=r-1=4k+2$,
and $X$ is contained in the cone with vertex a point over the $2$--uple
embedding of a threefold $Y$ as above (see Example \ref {exx}, (3));\pas

\item [(4)] $r=4k+3$, $\delta_k(X)=1$, $n_k=2$ and
$s^{(k)}(X)=r-1=4k+2$, and $X$ is the $2$--uple embedding
of a threefold $Y$ of degree $k$ in $\pp ^{k+1}$ with curve sections
of genus $0$, which is either a cone with vertex a line over a smooth
rational curve of degree $k$ in
$\pp^{k-1}$ or it has a double line (see Example \ref {exx}, (4)); \pas

\item [(5)]  $k=4$, $r=4k+3=19$, $\delta_4(X)=1$, $n_4=2$,
$s^{(k)}(X)=r-1=18$,  and $X$ is the $2$--uple embedding
of a threefold $Y$ in $\pp ^5$ with $\deg(Y)\geq 5$, contained in a
quadric (see Example \ref {exx}, (5));\pas

\item [(6)] $k\geq 4$, $r=4k+3$, $\delta_k(X)=1$, $n_k=2$,
$s^{(k)}(X)=r-1=4k+2$, and $X$ is the $2$--uple embedding of a threefold $Y$
of degree $k+1$ in $\pp^{k+1}$ with curve sections of arithmetic genus
2 (see Example \ref {exx}, (6));\pas 

\item [(7)] $r=4k+3-i$, $i=0,1$, $\delta_k(X)=1$,
$s^{(k)}(X)=r-1=4k+3-i$, $n_k=2-i$, and $X$ is contained in a cone with
vertex a space of dimension $k-i$ over the
$2$--uple embedding of a surface $Y$ of minimal degree $k$ in
$\pp^{k+1}$ (see Theorem \ref {dimY2}, case (1));\pas 

\item [(8)] $k=2$, $r=4k+3=11$, $\delta_k(X)=1$,
$s^{(k)}(X)=r-1=10$, $n_k=2$, and $X$ is contained in a cone with vertex
a line over the $2$--uple embedding of a
surface $Y$ of $\pp^3$ with $\deg(Y)\geq 3$ (see Theorem \ref {dimY2}, case
(2));\pas

\item [(9)] $k\geq 3$, $r=4k+3$, $\delta_k(X)=1$,
$s^{(k)}(X)=r-1=4k+2$, $n_k=2$, and  $X$ sits in a cone with
vertex of dimension
$k-1$ over the $2$--uple embedding of a surface $Y$ of degree $k+1$ in
$\pp^{k+1}$ with curve sections of arithmetic
genus $1$ (see Theorem \ref {dimY2}, case (3));\pas

\item [(10)]  $r\geq 4k+3$, $\delta_k=1$, $s^{(k)}(X)=4k+2$, $n_k=2$ and $X$ is
contained in a cone with vertex of dimension $k-1$, and not smaller, over a surface
which is not $k$--weakly defective (see Proposition \ref{proj});\pas

\item [(11)]  $r\geq 4k+3$, $\delta_k(X)=1$,
$s^{(k)}(X)=4k+2$, $n_k=2$ and $X$ is contained in a cone with vertex of
dimension $2k$, and not smaller, over a curve (see Proposition \ref{proj'});\pas

\item [(12)]  $r\geq 4k+2$, 
$s^{(k)}(X)=4k+1$, $n_k=1$ (hence $\delta_k(X)=1$, if $r=4k+2$, whereas 
$\delta_k(X)=2$, if $r>4k+2$) and
$X$ is contained in a cone with vertex of dimension $2k-1$, and not smaller, over a
curve (see Proposition \ref{proj''});\pas

\item [(13)] $4k+3\leq r\leq 4k+5$, $s^{(k)}(X)=4k+2$, $\delta_k(X)=1$,
$n_k=1$ and $X$ is either the $2$--uple
embedding of a threefold of minimal degree $k$ in $\pp^{k+2}$ (hence
$r=4k+5$), or the
projection from a point of $\pp^{4k+5}$ of the $2$--uple embedding of a threefold
$Y$ of minimal degree $k$ in $\pp^{k+2}$, or the projection from a line
$\ell\subset \pp^{4k+5}$ of the
$2$--uple embedding $Y'\subset \pp^{4k+5}$ of a threefold $Y$ of minimal
degree $k$ in
$\pp^{k+2}$ (see Example \ref {ex'});\pas

\item [(14)] $r=4k+3$, $s^{(k)}(X)=r-1=4k+2$, $\delta_k(X)=1$, $n_k(X)=2$ and $X$
is linearly normal, contained in the intersection of a
space of dimension $4k+3$ with the Segre embedding of
$\pp^{k+1}\times\pp^{k+1}$ in
$\pp^{k^2+4k+3}$, but not lying in the
$2$--uple embedding of
$\pp^{k+1}$, and such that the two projections of $X$ to $\pp^{k+1}$ span
$\pp^{k+1}$ (see Coorollary \ref {finalcor}).
\end{itemize}

Cases (1)--(9) correspond to the situation in which the tangential contact
locus is an
irreducible surface, case (10) corresponds to the situation in which
the tangential contact locus is a
reducible curve (a union of lines), cases (11)--(12) correspond to the
situation in which the tangential contact locus is a
reducible surface, cases (13)--(14) correspond to the situation in which the
tangential contact locus is
an irreducible curve, namely a rational normal curve of degree $2k$.\par

All threefolds in this list are actually minimally
$k$--defective. Threefolds of types (1)--(9) and (14) are not $h$--defective for
any $h>k$. The same is true for threfolds of type (13) with $r<4k+5$.
Threefolds of type (13) with $r=4k+5$ are also $(k+1)$--defective. Threefolds of
types (10)--(12) can be $h$--defective for $h>k$ (see Remark \ref {piudifett}).
\end{thm}

Brief discussions about the existence of smooth defective threefolds are contained
in Example \ref {exx} and Remarks \ref {liscezza} and \ref {piudifett} below.

\section  {PRELIMINARIES AND NOTATION} \label {notaz}

\subsection {} In this paper we work over the complex field $\mathbb C$. Let
$X\subseteq\pp^r$ be an irreducible projective scheme over $\C$. We will
denote by
$\deg(X)$ the {\it degree} of $X$ and  by $\dim(X)$ the {\it dimension} of
$X$. If $X$ is reducible, by $\dim(X)$ we mean the maximum of the
dimensions of its irreducible components. We will denote by $\mathcal
O_X$ the structure sheaf of $X$ and by $\mathcal I_X$ the ideal sheaf of
$X$ in $\p^r$.\par

If $X$ is irreducible, by a {\it general point} of $X$ we mean a point
which can vary in some dense open Zariski subset of $X$.

One says that $X$ is {\it rationally connected} is for $x,y\in X$ general points, there is a rational curve in $X$ containing them. \par 

If $Y\subset \pp^r$ is a subset, we denote by $<Y>$ the span of
$Y$. We will say that $Y$ is {\it non--degenerate} if
$<Y>=\p^r$. 

\subsection{} If $X\subseteq \p^r$ is an irreducible, projective, non
degenerate variety, we will say that it is a {\it cone} over a variety
$Y$, if there is a projective subspace $V$ of dimension
$v$, called a {\it vertex} of $X$, and a projective subspace $W$ of
dimension $r-v-1$, containing $Y$, such that $X=\cup_{x\in X,y\in
V}<x,y>$.\par

Let $V$ be a projective subspace of $\p^r$ of dimension
$v$, not containing $X$. 
We will say that {\it $X$ sits in the cone over a variety $Y\subset \p^{r-v-1}$ with vertex $V$} if and only if the image of the projection of
$X$ from $V$ is  $Y$.

\subsection{}\label{hilbfunc} Let $X\subset \p^r$ be an irreducible,
non--degenerate projective variety. One has
the following famous bound:

\begin{equation}\label{grado} \deg(X)\geq r-\dim(X)+1\end {equation}

\noindent which is sharp: the varieties achieving the bound \eqref
{grado} are called {\it varieties of minimal degree} and their
classification is well known  (see \cite {EH2}): curves of minimal degree
are the rational normal curves, surfaces of minimal degree are either
rational normal scrolls or the Veronese surface in $\p^5$, etc.\par

We will denote by $h_X$
the {\it Hilbert function} of $X$, namely for any non--negative integer
$k$, $h_X(k)$ is the dimension of the image of the restriction map:

$$\rho_{X,k}:H^0(\p^r, \mathcal O_{\p^r}(k))\to H^0(X, \mathcal O_X(k)).$$

Recall that one says that $X$ is $k$--normal if the map $\rho_{X,k}$ is
surjective, i.e. if $H^1(\p^r, \mathcal I_X(k))=0$. One says {\it linearly
normal}, {\it quadratically normal} etc. instead of $1$-normal,
$2$--normal, etc. Then $X$ is linearly normal if and only if the
hyperplane system $\mathcal H$ is complete. The variety $X$ is said to be
{\it projectively normal} if it is $k$--normal for every
$k\geq 1$.\par

We will need the following result of Castelnuovo's theory from \cite {WDV}
(see Theorem 6.1):

\begin{thm}\label {Hilb} Let $X\subset \p^r$ be an irreducible,
non--degenerate, projective variety of dimension $n$ and degree $d$. Set
$\iota:=\iota(X)=\min\{d+n-r-1,r-n\}$. Then:

\begin{equation}\label {cast} h_X(2)\geq \iota+r(n+1)-{\frac{n(n-1)}2}+1.\end{equation}

If $d\leq 2(r-n)+1$ then the following are equivalent:

\begin {itemize}

\item [(i)] the equality holds in \eqref{cast};
\item[(ii)] the  general curve section of $X$ is linearly normal of genus
$\iota$;
\item[(iii)] the  general curve section of $X$ is projectively normal of genus
$\iota$.
\end{itemize}

\end{thm}

Remark that, if $X$ is non--degenerate, the invariant $\iota(X)$  defined
in the previous statement is non--negative because of \eqref{grado}.

Let us now mention the following interesting proposition which follow from
results of
\cite {AR}, Corollary 3.3:

\begin{prop} \label {AlzRus} Let $X\subset \p^r$ be a smooth isomorphic projection of a
variety of minimal degree in $\p^{r+1}$. Then $X$ is $k$--normal for every
$k\geq 2$.\end{prop}

Using it, we can prove the following proposition, which
will be useful later:

\begin{lem}\label{RNS} Let $Y\subset \pp^{r}$ be a non--degenerate threefold
of minimal degree $r-2$ and let
$Y'\subset \pp^{r-1}$ be a projection of $Y$ from a point $p\notin Y$.
Assume that $Y'$ is non--singular in
codimension $1$. Then $h_{Y'}(2)\geq 4r-4$ and the equality holds if and
only if one of the following
cases occurs:

\begin{itemize}
\item [(i)] $Y'$ is a cone with vertex a line over a smooth rational curve
$C$ of degree $r-2$ in $\pp^{r-3}$;\pas
\item [(ii)] $Y'$ has a line $L$ of double points, and the pull back of
$L$ on $Y$ is a conic whose plane contains the centre
of projection $p$.\end{itemize}\end{lem}

\begin{proof} We let $S\subset \pp^{r-2}$ be the general surface section
of $Y'$ and
$C\subset \pp^{r-3}$ the general curve section of $S$. Notice that
$C$ is a smooth
rational curve of degree $r-2$. Standard diagram chasing gives:

$$h_{Y'}(2)-h_{Y'}(1)=h_S(2)+\dim(\ker\{H^1(\p^{r-1},{\mathcal
I}_{Y'}(1))\to H^1(\p^{r-1}, {\mathcal I}_{Y'}(2))\})$$
$$h_{S}(2)-h_{S}(1)=h_C(2)+\dim(\ker\{H^1(\p^{r-2},{\mathcal
I}_{S}(1))\to H^1(\p^{r-2},{\mathcal I}_{S}(2)) \})$$

\noindent (see \cite {Ci}, p. 30). Now:

$$h_{Y'}(1)=r, \quad h_{S}(1)=r-1$$
and

$$h_C(2)=2r-3$$
by Proposition \ref {AlzRus}. Thus:

$$ h_{Y'}(2)=4r-4+$$

\begin{equation} \label{delta'}
+\dim(\ker\{\p^{r-1}, H^1({\mathcal I}_{Y'}(1))\to H^1(\p^{r-1},{\mathcal
I}_{Y'}(2))\})+\end{equation}

$$+\dim(\ker\{H^1(\p^{r-2},{\mathcal I}_{S}(1))\to
H^1(\p^{r-2},{\mathcal I}_{S}(2))\}).$$

This proves the first assertion. Moreover $h_{Y'}(2)=4r-4$ yields:

$$\dim(\ker\{H^1(\p^{r-1},{\mathcal I}_{Y'}(1))\to H^1(\p^{r-1},{\mathcal
I}_{Y'}(2))\})=0$$

$$\dim(\ker\{\p^{r-2},H^1({\mathcal I}_{S}(1))\to
H^1(\p^{r-2},{\mathcal I}_{S}(2)) \})=0.$$

This implies that both $Y'$ and $S$ are singular. Suppose in fact that
$Y'$ is smooth. Then
$h^1(\p^{r-1},{\mathcal I}_{Y'}(1))=1$ and $h^1(\p^{r-1},{\mathcal
I}_{Y'}(2))=0$ by Proposition \ref {AlzRus}. The
same argument works for $S$. \par

Let $x\in S$ be a singular point. If $S$ is a cone with vertex $x$ we
are in case (i). Suppose $S$ is
not a cone and let $\mu$ be the multiplicity of $S$ at $x$. The
projection of $S$ from $x$ is
a non--degenerate surface of degree $r-2-\mu$ in $\pp^{r-3}$. This
proves that $\mu=2$.
We claim that $x$ is the only
singular point of $S$. In fact $S$ is the
projection of a surface $T\subset \pp^{r-1}$ of minimal degree from
a point $p\notin T$. Since $S$ is not a cone, then $T$ is not a cone, hence
it is smooth (see \cite {EH2}). The
singular point $x$ of $S$ arises from a secant (or tangent) line $\ell$ to $T$ passing
through
$p$. Suppose there is another singular
point $y\in S$. This would correspond to some other secant line $\ell'$
to $T$  containing $p$. The
plane $\Pi=<\ell,\ell'>$ would then be $4$--secant to $T$, and this is
possibe only if $\Pi$ intersect
$T$ along a conic (see again \cite {EH2}). This would in turn yield a
double line for $S$ and codimension $1$ singularities
for $Y$, a contradiction. In conclusion we are in case (ii).\par

Conversely, if we are in case (i), one has

$$h^0(\p^{r-1},{\mathcal I}_{Y'}(2))=h^0(\p^{r-3},{\mathcal I}_C(2))={\binom{r-1}
2}-2r+3$$ (see again Proposition \ref {AlzRus}). Hence:

$$h_{Y'}(2)={\binom{r+1} 2}-{\binom{r-1} 2}+2r-3=4r-4.$$

If we are in case (ii), we have $h^1(\p^{r-1},{\mathcal
I}_{Y'}(1))=h^1(\p^{r-2},{\mathcal I}_{S}(1))=0$, because the
singular varieties $Y'$ and $S$ are linearly normal, and formula \eqref
{delta'} yields
$h_{Y'}(2)=4r-4$. \end{proof}

It is clear that one can prove a similar more general result for
projections of varieties of minimal degree of any dimensions. Since we will
not need it, we do not dwell on this here.\par

We also record the following:

\begin{prop} \label {ratsing} Let $Y\subset \p^{k+1}$ be an irreducible,
non--degenerate, threefold of degree $k$ with smooth curve sections of
genus $0$ and a singular line. Then:

\begin{itemize}

\item[(i)] either $Y$ is a cone over a smooth rational curve of degree $k$
in $\p^{k-1}$,

\item [(ii)] or $Y$ is the  projection in $\pp^{k+1}$ of a
threefold $Z$ of minimal degree $k$ in $\pp^{k+2}$ from a point $p\notin
Z$, containing a conic sitting in a plane passing through $p$.
\end{itemize}

\end{prop}

\begin{proof} One knows that the threefold $Y$ is the  projection in $\pp^{k+1}$ of a
threefold $Z$ of minimal degree $k$ in $\pp^{k+2}$ from a point $p\notin
Z$. If $Z$ is a cone with vertex a line, we are in case (i). If $Z$ is a
cone with vertex a point, then $Y$ is also a cone. However the point $p$
has to sit on a secant line $\ell$ to $Z$ and we are in case (ii): the conic
here is reducible in the two generators of the cone $Z$ passing through the
intersection points of $\ell$ with $Z$. If $Z$ is smooth, then again $p$
has to sit on a secant line $\ell$ to $Z$. The argument to show that we are
in case (ii) then goes as in the proof of  Lemma \ref {RNS}.
\end{proof}

\subsection{}\label {linsyst} Let $X$ be an irreducible, projective
variety. We will use the symbol $\equiv$ to denote linear equivalence 
of Weil divisors, or linear systems, on $X$. If $D$ is a divisor, we will denote,
as usual, by $|D|$ the complete linear system of $D$. \par

Let $H$ be a hyperplane divisor on $X$. We will denote by ${\mathcal
H}\subset |H|$ the
(possibly not complete) {\it hyperplane system}, i.e. the linear system cut
out on
$X$ by the hyperplanes of $\pp^r$. If
we  assume that $X$ is non--degenerate, then ${\mathcal H}$ has dimension
$r$. In this case, by abusing notation, we will sometimes identify the
divisor $H$ with the unique hyperplane which cuts $H$ on $X$.\par

If $\mathcal L$ is a linear system of dimension $r$ of Weil
divisors on $X$, we will denote by:

$$\phi_\mathcal L: X \map \pp^{r}$$

\noindent the rational map defined by $\mathcal L$. \par

If $p_1,...,p_k\in X$ are smooth points and $m_1,...,m_k$ are positive integers,
we will denote by 
$\mathcal L(-m_1p_1-...-m_kp_k)$ the sublinear system of $\mathcal L$ formed by
all divisors in
$\mathcal L$ with multiplicity at least $m_i$ at $p_i$, $i=1,...,k$.\par

If ${\mathcal L}_1$ and ${\mathcal L}_2$ are linear
systems of Weil divisors
on a $X$, we define ${\mathcal L}_1+{\mathcal L}_2$ as the
minimal linear system
of Weil divisors on $X$ containing all divisors of the form $D_1+D_2$
where $D_i\in {\mathcal
L}_i$, $i=1,2$.  The linear system ${\mathcal L}_1+{\mathcal L}_2$ is
called the {\it
minimal sum} of ${\mathcal L}_1$ and ${\mathcal L}_2$. \par

If ${\mathcal L}$ is a linear system on $X$, one writes $2{\mathcal L}$
instead of
${\mathcal L}+{\mathcal L}$. Similarly one can consider the linear
system  $h{\mathcal
L}$ for all positive integers $h$.\par

Let ${\mathcal L}_1$ and ${\mathcal L}_2$ be linear
systems of Weil divisors
on $X$ of dimensions $r_1, r_2$.  Set ${\mathcal L}={\mathcal L}_1+{\mathcal
L}_2$, and suppose $\mathcal L$ has dimension $r$. One can consider the
maps:

$$\phi_{{\mathcal L}_i}: V\to \pp^{r_i}, i=1,2$$

$$\phi_{{\mathcal L}}: V\to \pp^{r}$$

\noindent determined by the linear systems in question. It is clear that:

$$\phi_{{\mathcal L}}= \psi\circ (\phi_{{\mathcal L}_1}\times
\phi_{{\mathcal L}_2})$$
where:

$$\psi: \pp^{r_1}\times \pp^{r_2}\to \pp^{r_1r_2+r_1+r_2}$$
is the {\it Segre embedding}.

Similarly if ${\mathcal L}$ is a linear system on $V$ with dimension $r$,
and

$$\phi_{{\mathcal L}}: V\to \pp^{r}$$
is the corresponding map, then for any positive integer $h$ one
has:

$$\phi_{h{\mathcal L}}=\psi_h\circ \phi_{{\mathcal L}}$$
where:

$$\psi_h: \pp^{r}\to \pp^{\binom{r+m} m -1}$$
is the {\it $h$-th Veronese embedding}. 

We will need the following lemma:

\begin{lem}\label{hopf} Let $C$ be a smooth, irreducible, projective
curve of genus $g$ and let
$\mathcal L_i$, $i=1,2$, be base point free linear systems on $C$. Set
$\mathcal L=\mathcal L_1+\mathcal
L_2$, $\dim(\mathcal L)=r$, $\dim(\mathcal L_i)=r_i$, $i=1,2$. One has $r\geq
r_1+r_2$. If $\mathcal L$ is birational and $r=r_1+r_2$, then $g=0$.
If moreover $r_1=r_2$, then $\mathcal L_1=\mathcal
L_2$ and this series is complete.\end{lem}

\begin{proof} The inequality $r\geq r_1+r_2$ follows from Hopf's lemma
(see \cite {ACGH}, p. 108). If
the equality holds, then $g=0$ by \cite {ACGH}, Exercise B-1, p. 137.
Suppose now that, in addition, $r_1=r_2$.
Then by \cite {ACGH}, Exercise B-4, p. 138, the two series are equal, and by
the birationality assumption on $\mathcal L$, also $\mathcal L_1=\mathcal
L_2$  is birational. If the series is not complete, then its image is a
non-normal
rational curve $\Gamma$ in $\pp^{r_1}$. Then, by Theorem \ref {Hilb},
$h_\Gamma(2)\geq 2r_1+2$, a contradiction.  \end{proof}

Let us now recall a definition:

\begin{defn}\label{bigdef}\rm We will say that a linear system $\mathcal L$
on a smooth, projective
variety $Z$ is {\it birational} [resp. {\it very big}] if the map
$\phi_{\mathcal L}$ determined by
$\mathcal L$ is birational [resp. generically finite] to its image. The
system $\mathcal L$ is {\it
big} if the complete system determined by some multiple of $\mathcal L$
is big, or, equivalently,
if the complete system determined by some multiple of $\mathcal L$ is
birational. \end{defn}

The following definition is standard:

\begin{defn}\rm If $W$, $Z$ are varieties with $W\subset Z$ and
$\mathcal L$ is a linear system of
Cartier divisors on $Z$, one denotes by $\mathcal L_{|W}$ the {\it
restriction} of $\mathcal
L$ to $W$ defined as follows. The system $\mathcal L$  corresponds to a
vector subspace $V$ of
$H^0(Z,\mathcal O_Z(L))$, with $L$ a divisor in $\mathcal L$, and
$\mathcal L_{|W}$ corresponds to the
image of  $V$ via the restriction map $H^0(Z,\mathcal O_Z(L))\to
H^0(W,\mathcal O_W(L))$.\end{defn}

The following lemma is immediate and the proof can be left to the reader:\par

\begin{lem}\label{big} Let $Z$ be a smooth, projective variety, let
$\mathcal L$ be a linear system on
$Z$ and let $\mathcal V$ be a family of closed subvarieties of $Z$ such
that a general element $V$ of
$\mathcal V$ is smooth and irreducible. Suppose that if $z\in Z$ is a
general point, there is a
variety $V$ in $\mathcal V$ containing $z$. Then if $\mathcal L$ is
birational [resp. very
big, big] and $V$ is a
general element of
$\mathcal V$, then $\mathcal L_{|V}$ is also birational [resp. very big,
big].\end{lem}

\subsection {} \label {monodromy}  Let $X$ be an
irreducible projective variety. For any positive integer $h$ we let
${\rm Sym}^h(X)$ be the $k$--fold symmetric product of $X$. If
$p_1,...,p_h$ are points in $X$, we denote by $[p_1,...,p_h]$ the
corresponding point in ${\rm Sym}^h(X)$. Namely one has a
surjective morphism:

$$\pi_{X,h}: (p_1,...,p_h)\in X^h\to [p_1,...,p_h]\in {\rm Sym}^h(X)$$

\noindent which is a finite covering of degree $h!$ and the monodromy, or
Galois, group of this covering is the full symmetric group $\mathcal S_h$.

\subsection{}\label{invo} Let $X$ be an irreducible, projective variety.
Let $\mathcal D=\{D_y\}_{y\in Y}$  be an algebraic family
of Weil divisors on $X$ parametrized by an projective variety $Y$. We will
constantly assume that
$\mathcal D$ is {\it effectively parametrized} by $Y$, i.e. that the
corresponding map of $Y$ to the appropriate Hilbert scheme of divisors on
$X$ is generically finite. One says that $\mathcal D$ is {\it irreducible
of dimension $m$} if $Y$ is.\par

An irreducible algebraic family $\mathcal D$ of dimension $m$ on $X$ is
called an {\it involution} if there is one single divisor of $\mathcal D$
coontaining $m$ general points of $X$.\par 

Involutions have been studied in \S 5 of \cite{WDV} to which we defer the
reader for details. In particular in \cite{WDV} a classical theorem of
Castelnuovo and Humbert concerning involutions on curves (see
\cite{WDV}, Proposition 5.9),
has been extended to higher dimensional varieties (see
\cite{WDV}, Proposition 5.10). For the reader's
convenience we recall here the result we will need later on:

\begin{thm}\label{CasHum} Let $X$ be an irreducible, projective variety of
dimension $n>1$. Let $\mathcal
D$ be an $m$--dimensional involution with no fixed divisors and such that
its general element
is reduced. Then either $\mathcal D$ is a linear system of Weil divisors
or it is {\rm composed
with a pencil}, i.e. there is a map $f: X\map C$ of $X$ to an irreducible
curve $C$ such that
$\mathcal D$ is  the
pull-back, via
$f$, of an involution on $C$. If $m>1$, the general element of $\mathcal D$
is reducible if and only if $\mathcal D$ is composed with a pencil.
\end{thm}

\subsection {}\label{subsecsec} Let $X\subset \pp^r$ an irreducible, non--degenerate projective variety of dimension $n$. Let $k$ be a non--negative integer and
let $S^{k}(X)$ be the $k$--{\it secant variety} of $X$, i.e. the
Zariski closure in $\pp^r$ of the set:\medskip
 
\centerline {$\{ x\in \pp^r: x$  lies in the span of $k+1$
independent points of $X\}$}\medskip
 
Of course $S^0(X)=X$, $S^r(X)=\pp^r$ and $S^{k}(X)$ is empty if
$k\geq r+1$. Moreover, for any $k\geq 0$, one has $S^{k}(X)\subseteq
S^{k+1}(X)$. We will write
$S(X)$ instead of
$S^1(X)$ and we will assume $k\leq r$ from now on. \par
 
One can consider the {\it abstract $k$--th secant variety} $S^k_X$
of $X$, i.e. $S^{k}_X\subseteq {\rm Sym}^k(X)\times\pp^r$ is the
Zariski closure of the set of all pairs $([p_0,...,p_k],x)$ such
that $p_0,...,p_k\in X$ are linearly independent points and $x\in
<p_0,...,p_k>$. One has the surjective map $p_X^k:S^k_X\to
S^k(X)\subseteq\pp^r$, i.e. the projection to the second factor.
Hence:
 
\begin{equation} \label{defect} s^{(k)}(X):= \dim (S^{k}(X))\leq
\min\{r,\dim(S^k_X)\}= \min\{r,n(k+1)+k\}\end{equation}
 
The right hand side of \eqref {defect} is called the {\it expected
dimension} of $S^{k}(X)$ and will be denoted by $\sigma^{(k)}(X)$.
One says that $X$  is $k$-{\it defective} when strict
inequality holds in (\ref {defect}). One says that:
 
$$\delta_k(X):=\sigma^{(k)}(X)-s^{(k)}(X)$$
 
\noindent is the $k$--{\it defect} of $X$. The variety $X$ is called 
{\it defective} if it is $k$--defective
for some $k\geq 1$. Observe that, if $r\geq n(k+1)+k$, then $X$ is $k$--defective if and only if there are infinitely many $(k+1)$--secant $\p^ k$'s passing through the general point of $S^ k(X)$.
Notice also that, if $X$ is $h$--defective for some $h\geq
1$, then it
is $k$--defective for all $k$ such that $k\geq h$ and $s^{(k)}<r$. If $X$ is
$k$--defective but not $(k-1)$--defective, then we will say that $X$ is
{\it minimally $k$--defective}. \par

We will write $s^{(k)}, \sigma^{(k)}, \delta_k$ etc. instead of $s^{(k)}(X),
\sigma^{(k)}(X), \delta_k(X)$, if there is no danger of confusion.

\subsection {}\label {terr} If $p$ is a smooth point of $X$, we denote
by $T_{X,p}$ the {\it tangent
space} to $X$ at $p$. If $\Pi$ is a projective subspace of $\p^r$, we
say that $\Pi$ is {\it tangent} to $X$ at $p$ if either: 

\begin{itemize}

\item  $\dim(\Pi)\leq \dim(X)$ and  $p\in \Pi\subseteq T_{X,p}$, or

\item  $\dim(\Pi)\geq \dim(X)$ and $T_{X,p}\subseteq \Pi$.
\end{itemize}

Let $k$ be a positive integer and
let $p_1,...,p_k$ be points of $X$. We denote by
$T_{X,p_1,...,p_k}$ the span of $T_{X,p_i}, i=1,...,k$.\par
 
If $X\subset\pp^r$ is a projective variety, Terracini's Lemma
describes the tangent space to $S^k(X)$ at a general point of it and
gives interesting information in case $X$ is $k$--defective
(see \cite {Terr1} or, for modern versions, \cite {Adl}, \cite
{WDV}, \cite {Dale1}, \cite {Zak}). We may state it as follows:
 
\begin{thm}\label{terracini} (Terracini's Lemma) Let $X\subset\pp^r$
be an irreducible, projective variety. If $p_0,...,p_k\in X$ are
general points and $x\in <p_0,...,p_k>$ is a general point, then:
 
$$T_{S^k(X),x}=T_{X,p_0,...,p_k}.$$
 
If $X$ is $k$--defective, then:
\begin{itemize}

\item [(i)]  $T_{X,p_0,...,p_k}$ is tangent to $X$ along a variety
$\Gamma:=\Gamma_{p_0,...,p_k}$ of positive dimension
$\gamma_k:=\gamma_k(X)$  containing $p_0,...,p_k$; 

 \item [(ii)] the general hyperplane $H$
containing $T_{X,p_0,...,p_k}$ is tangent to $X$ along a variety
$\Sigma:=\Sigma(H):=\Sigma_{p_0,...,p_k}(H)$ of  positive dimension
$\epsilon_k:=\epsilon_k(X)$ containing $p_0,...,p_k$. 
\end{itemize}

One has $\Gamma\subseteq \Sigma$  and
therefore $\gamma_k\leq \epsilon_k$. Moreover one has:
 
\begin {equation}\label{terracini'} k\leq \dim (<\Gamma>)\leq
\dim(<\Sigma>)\leq k\epsilon_k+k+\epsilon_k-\delta_k.\end{equation}
 
\end{thm}

We may and will assume that all irreducible components of
$\Gamma$ and of $\Sigma$ contain some of the points
$p_0,...,p_k$. Otherwise we simply get rid of those components that do not
do so. We will call $\Gamma$ the {\it tangential
$k$--contact locus} of $X$ at $p_0,...,p_k$. Similarly, for any hyperplane
$H$ containing $T_{X,p_0,...,p_k}$, we will call $\Sigma:=\Sigma(H)$ the 
{\it $k$--contact locus} of $H$ at $p_0,...,p_k$. We will call
$\epsilon_k(X)$ the $k$-{\it singular defect} of  $X$ and $\gamma_k(X)$
the $k$-{\it tangential defect} of $X$. 

The following is a well known, straightforward application of Terracini's
lemma (see \cite {Zak}):

\begin{prop}\label{crescente} Let $X\subset \p^r$ be an irreducible,
projective, non--degenerate variety. If for some $k\geq 0$ one has 
$s^{(k)}(X)=s^{(k+1)}(X)$, then $s^{(k)}(X)=r$. In particular, if
$s^{(k)}(X)=r-1$, then $s^{(k+1)}(X)=r$.

\end{prop}

We record also the following result whose easy proof can be left to the
reader:

\begin{prop}\label {proiezione} Let $X\subset \p^r$ be an irreducible,
projective, non--degenerate variety of dimension $n$. let $\Pi\subset \p^r$
be a projetive subspace of dimension $s$ and let $\pi$ be the projection of
$\p^r$ from $\Pi$ to $\p^{r-s-1}$. Let $Y=\pi(X)$ and let $m$ be its
dimension. Then:

\begin {itemize}

\item[(i)] for every positive integer $k$ one has $\pi(S^k(X))=S^k(Y)$,
hence $s^{(k)}(Y)\leq s^{(k)}(X)\leq s^{(k)}(Y)+s+1$;

\item [(ii)] if $n=m$ and $s^{(k)}(Y)=(k+1)n+k$ then
also $s^{(k)}(X)=(k+1)n+k$;

\item [(iii)] if $n=m$, $s^{(k)}(Y)=s^{(k)}(X)$ and $Y$ is $k$--defective
[resp. minimally $k$--defective] then also $X$ is $k$--defective
[resp. minimally $k$--defective] and $\pi$ maps the tangential $k$--contact
locus [resp. the $k$--contact locus] of $X$ to the tangential $k$--contact
locus [resp. the $k$--contact locus] of $Y$.
\end{itemize}

\end{prop}

\subsection{}\label{subsecwd} We recall from \cite {WDV} the definition of
a {\it $k$--weakly defective} variety, i.e. a variety
$X\subset\p^r$ such that if $p_0,...,p_k\in X$ are general points,
then  the general hyperplane $H$ containing $T_{X,p_0,...,p_k}$ is
tangent to $X$ along a variety $\Sigma:=\Sigma(H):=\Sigma_{p_0,...,p_k}(H)$
of positive dimension $\epsilon_k:=\epsilon_k(X)$ containing $p_0,...,p_k$.
Similarly we  can say that a variety
$X\subset\p^r$ is {\it $k$--weakly tangentially defective} if whenever
$p_0,...,p_k\in X$ are general points, then $T_{X,p_0,...,p_k}$ is
tangent to $X$ along a variety $\Gamma:=\Gamma_{p_0,...,p_k}$ of 
positive dimension $\gamma_k:\gamma_k(X)$ containing $p_0,...,p_k$. Of
course 
$k$--weakly tangentially defectiveness implies $k$--weakly
defectiveness, but the converse does not hold in general. Moreover, by
Terracini's lemma, a $k$--defective variety is also $k$--weakly
defective but again the converse does not hold in general (see \cite
{WDV}).
 
\begin{rem} \label{0wd} {\rm  A curve which is not a line is never
$k$--weakly (tangentially) defective for any $k$. Hence a curve is
never $k$--defective.

For a variety $X$ being
$0$--tangentially defective means that it is
{\it developable}, i.e. the Gauss map of
$X$ has positive dimensional fibres, in particular, according to Zak's
 theorem on tangencies (see \cite  {Zak}),
$X$ must be singular, unless $X$ is a linear space. Instead, $0$--weakly
defective means that the dual variety of $X$ is not a hypersurface. 
In the surface case this happens if and only if the surface is
developable, i.e. if and only if the surface is either a
cone or the tangent developable to a curve (see \cite {GrHarr}). However
this is no longer the case if $\dim(X)>2$ (see \cite {ein}). Finally no
variety is
$0$--defective.}\end{rem}

\subsection {}\label {tangentialproj} Let $X\subset \p^r$ be, as above,
an irreducible, non--degenerate, projective variety of dimenszion $n$.  For
every non--negative integer $k$, we set:

$$r_k=r_k(X):=r-\dim(T_{X,p_1,\dots,p_k})-1=r-s^{(k-1)}(X)-1.$$

Consider the projection of $X$
with centre
$T_{X,p_1,...,p_k}$. We call this a {\it general $k$--tangential
projection} of $X$, and we will denote it by $\tau_{X,p_1,...,p_k}$ or
simply by $\tau_{X,k}$. We will denote by $X_k$ its image. Notice that $X_k$ is
non--degenerate in $\pp^{r_k}$. We define $\tau_{X,0}$ as the
identity so that $X_0=X$.\par

We set $n_k:=n_k(X):=\dim(X_k)$ and $m_k:=m_k(X)=n-n_k$. Notice that $m_k$
is the dimension of the general fibre of the map $\tau_{X,k}: X\map
X_k$.\par

\begin{lem}\label {emme}  Let $X\subset \p^r$ be 
an irreducible, projective variety of
dimension $n\geq 2$. Then $m_k\leq \gamma_k$.
\end{lem}

\begin{proof} Consider the
general $k$--tangential
projection $\tau_{X,k}: X\map X_k$ from
$T_{X,p_1,...,p_h}$.
Let $p_0$ be a general point of $X$. The pull--back via
$\tau_{X,k}$ of the tangent space to $X_k$
at $\tau_{X,h}(p_0)$ is $T_{X,p_0,\dots,p_h}$.
Hence $T_{X,p_0,\dots,p_k}$ is tangent to
$X$ along the whole fibre $\tau_{X,k}^{-1}(\tau_{X,k}(p_0))$, which has
dimension $m_k$. By the definition of $\gamma_k$ we have the
assertion.\end{proof}

The following result is an immediate consequence of Terracini's lemma
(see \cite {WDV}, \S 3):

\begin{prop}\label {difetti} Let $X\subset \p^r$ be an irreducible,
projective variety of dimension $n\geq 2$ which is minimally
$k$--defective. Then:

\begin{itemize}

\item [(i)]
$n_h=n$ for 
$1\leq h\leq k-1$, whereas
$n_k\leq n-\delta_k<n$, thus $m_h=0$ for 
$1\leq h\leq k-1$, whereas $m_k\geq \delta_k$;

\item [(ii)] $0<n_k<r_k$, i.e. $X_k$ is a proper subvariety of positive
dimension of
$\p^{r_k}$;

\item [(iii)] if $r\geq n(k+1)+k$ then
$n_k=n-\delta_k$, i.e.
$m_k=\delta_k$;

\item [(iv)] $\delta_k\leq n-1$ and
$r\geq (n+1)k+2$;

\item [(v)] if $r=(n+1)k+2$ then
$\delta_k=1$, $m_k=n-1$ and $X_k$ is a plane curve.

\end{itemize}
\end{prop}

\begin{proof} Let us prove part (i). 
Consider the
general $h$--tangential
projection $\tau_{X,h}: X\map X_h\subseteq \pp^{r_h}$ from
$T_{X,p_1,...,p_h}$.
Let $p_0$ be a general point of $X$. For all $h$, the pull--back via
$\tau_{X,h}$ of the tangent space to $X_h$
at $\tau_{X,h}(p_0)$ is $T_{X,p_0,\dots,p_h}$, hence:

\begin{equation}\label{tangsp}
s^{(h)}=\dim(T_{X,p_0,\dots,p_h})=n_h+\dim(T_{X,p_1,\dots,
p_h})+1=n_h+s^{(h-1)}+1.
\end{equation}

Since $X$ is minimally $k$--defective, then for all $i<k$ one has
$s^{(i)}= (i+1)n+i$.
Hence formula (\ref {tangsp}) gives, for all $h<k$, 
$n_h=n$. \pas

On the other hand, we know that:

\begin{equation}\label{tangsp'}
\dim(T_{X,p_0,\dots,p_k})=s^{(k)}<\sigma^{(k)}\leq n(k+1)+k.
\end{equation}
 
Hence by formula  (\ref {tangsp}), applied for $h=k$, we have that
$n_k\leq n-\delta_k<n$. This proves (i).\par

Since $X$ is $k$-defective, then
$\dim(S^k(X))<r$, hence $\dim(T_{X,p_0,\dots,p_k})=s^{(k)}<r$, i.e.
$T_{X,p_0,\dots,p_k}$ cannot
coincide with the whole space. Therefore $X_k$ is a proper subvariety of
$\pp^{r_k}$. Since it is non--degenerate, one has $n_k>0$. This proves part
(ii).\par

If $r\geq n(k+1)+k$ then $\sigma^{(k)}= n(k+1)+k$. Therefore
\eqref{tangsp} and \eqref{tangsp'} yield $n_k= n-\delta_k$, i.e.
part (iii).

Since 
$n_k>0$ one has
$\delta_k<n$. Furthermore $r-(n+1)k=r_k>n_k\geq 1$. Hence $r\geq
(n+1)k+2$. This proves part (iv).

If $r=(n+1)k+2$, then
$\sigma^{(k)}=\min\{n(k+1)+k,r\}=r=(n+1)k+2$ and 
formula  (\ref {tangsp}), applied for $h=k$, gives
$n_k=2-\delta_k$. On the other hand
 $2=r-(n+1)k=r_k>n_k\geq 1$, thus $n_k=1$ and
$\delta_k=1$. This proves part (v).\end{proof}

The following example is well known:

\begin{ex}\label{segre}\rm  Let $k\geq 1$ be an integer and let $X$ be the Segre
embedding of $\p^{k+1}\times
\p^{k+1}$ in $\p^{k^2+4k+3}$. We claim that $X$ is $1$--defective with
$\delta_1(X)=2$, $n_1(X)=2k$.\par

Indeed one sees that $X_1$ is nothing but the Segre
embedding of $\p^{k}\times
\p^{k}$ in $\p^{k^2+2k}$. The assertion immediately follows. 
\end{ex}

\subsection {} \label {incones} We recall the following elementary
criterion, needed in the sequel, which tells us when a variety sits in a
cone (see \cite {WDV}, Proposition 4.1):

\begin{prop} \label {varincon} Let $X\subset \p^r$ be an irreducible,
projective variety of dimension $n$. Let $\Pi\subset \p^r$ be a projective
subspace of dimension $s$ not containing $X$. Then $X$ projects from $\Pi$
to a variety $Y$ of dimension $m<n$, i.e. $X$ sits in the cone with vertex
$\Pi$ over $Y$, if and only if the general tangent space to $X$ intersects
$\Pi$ along a subspace of dimension $n-m-1$. In particular $X$ is a
cone with vertex $\Pi$ if and only if the general tangent space to $X$
contains $\Pi$ and $X$ sits in a $\p^{s+1}$ containing 
$\Pi$ if and only if the general tangent space to $X$ meets $\Pi$ along a 
a subspace of dimension $n-1$.\end {prop}

As a consequence we have:

\begin{prop}\label {conidif}
Let $X\subset \p^r$ be an irreducible,
projective variety of dimension $n$. Let $\Pi\subset \p^r$ be a projective
subspace of dimension $s\geq 0$ not containing $X$, let $\Pi'$ be a complementary subspace of dimension $r-s-1$, and assume that $X$
projects from $\Pi$ to a variety $Y\subset \Pi'$ of dimension  $n-1$.
Then:

\begin{itemize}

\item [(i)] for every $k\geq s$, $S^k(X)$ is the cone with vertex
$\Pi$ over
$S^k(Y)$, and therefore $s^{(k)}(X)=s^{(k)}(Y)+s+1$ and 
$\gamma_k(X)=\gamma_k(Y)+1$. In particular $X$
is $k$ defective if $k>s$ [and $\gamma_k(X)=1$ if and only if $Y$ is not 
$k$--weakly defective] and $X$ is $s$--defective if and only if $Y$ is
$s$--defective;

\item [(ii)] for
every positive integer $k<s$, one has $s^{(k)}(X)\geq s^{(k)}(Y)+k+1$. In
particular if $Y$ is not $k$--defective, then the equality holds and
$X$ is also not $k$--defective. 
\end{itemize}

As a consequence, if $Y$ is not $s$--defective, then $X$ is minimally
$(s+1)$--defective, whereas if $Y$ is minimally
$s$--defective, then also $X$ is minimally $s$--defective. Finally then $\gamma_k(X)=1$ if and only if $Y$ is not 
$(s+1)$--weakly defective. 

\end{prop} 

\begin{proof} Let $p\in X$ be a general point. By Proposition \ref
{varincon}, $T_{X,p}$ meets $\Pi$ in one point. Hence there is a map:

$$f: p\in X\map  T_{X,p}\cap \Pi\in \Pi$$

\noindent and the image of $f$ is non--degenerate in $\Pi$ since $X$ is
non--degenerate in $\p^r$. \par

Let $k\geq s$. Let $p_0,...,p_k$ be general points of $X$ and let $q_0,...,q_k$ be the corresponding projections on $Y$. Then $T_{X,p_0,...,p_k}$
contains $\Pi$ and projects from $\Pi$ to $T_{Y,q_0,...,q_k}$. 
Hence the general tangential $k$--contact locus of $S$  pulls--back to $X$ to the general tangential $k$--contact locus of $X$. Moreover, by Terracini's lemma and by Proposition \ref {varincon}, we
have that $S^{(k)}(X)$ is a cone with vertex $\Pi$. By Proposition \ref
{proiezione}, $S^{(k)}(X)$ is the cone with vertex $\Pi$ over
$S^{(k)}(Y)$. Thus (i) follows.\par

Let $k\leq s$ be a positive integer. The above argument shows that
$T_{X,p_0,...,p_k}$ intersects $\Pi$ along a $\p^{l}$, with $l\geq k$.
Since $S^{k}(X)$ projects from $\Pi$ to
$S^{k}(Y)$, part (ii) follows from Proposition \ref {varincon}.\par

The rest of the assertion is trivial. \end{proof}

In a similar way one proves the following:

\begin{prop}\label {conidif'}
Let $X\subset \p^r$ be an irreducible,
projective variety of dimension $n$. Let $\Pi\subset \p^r$ be a projective
subspace of dimension $2s\geq 0$ not containing $X$ and assume that $X$
projects from $\Pi$ to a variety $Y\subset \p^{r-2s-1}$ of dimension  $n-2$.
Then:

\begin{itemize}

\item [(i)] for every $k\geq s$, $S^{(k)}(X)$ is the cone with vertex
$\Pi$ over
$S^{(k)}(Y)$, hence $s^{(k)}(X)=s^{(k)}(Y)+2s+1$ and therefore $X$ is $k$ defective;

\item [(ii)] for
every positive integer $k<s$, one has $s^{(k)}(X)\geq s^{(k)}(Y)+2k+2$. In
particular if $Y$ is not $k$--defective, then the equality holds and
$X$ is also not $k$--defective. 
\end{itemize}

As a consequence, if $Y$ is not $(s-1)$--defective, then $X$ is minimally
$s$--defective. 
\end{prop} 

\begin{prop}\label {conidif''}
Let $X\subset \p^r$ be an irreducible,
projective variety of dimension $n$. Let $\Pi\subset \p^r$ be a projective
subspace of dimension $2s-1\geq 0$ not containing $X$ and assume that $X$
projects from $\Pi$ to a variety $Y\subset \p^{r-2s}$ of dimension  $n-2$.
Then:

\begin{itemize}

\item [(i)] for every $k\geq s-1$, $S^{(k)}(X)$ is the cone with vertex
$\Pi$ over
$S^{(k)}(Y)$, hence $s^{(k)}(X)=s^{(k)}(Y)+2s$.  In particular $X$
is $k$ defective if $k>s-1$ and $X$ is $(s-1)$--defective if and only if $Y$ is
$(s-1)$--defective;

\item [(ii)] for
every positive integer $k<s-1$, one has $s^{(k)}(X)\geq s^{(k)}(Y)+2k+2$. In
particular if $Y$ is not $k$--defective, then the equality holds and
$X$ is also not $k$--defective. 
\end{itemize}

As a consequence, if $Y$ is not $(s-1)$--defective, then $X$ is minimally
$s$--defective. 
\end{prop} 

\subsection{} Now we give a definition which extends a
definition given in \cite{Sacch}.

\begin{defn} \rm Let $Z\subseteq \pp^r$ be an irreducible,
non--degenerate, projective variety of
dimension $n$. Let $h$ be a non--negative integer such that $h\leq n$.
Let $z\in Z$ be a general point
and suppose that $\overline {T_{Z,z}\cap Z-\{x\}}$ has an irreducible
component of dimension $h$. In
this case we will say that $Z$ is {\it $h$--tangentially degenerate}. We
will simply say that $Z$ is
{\it tangentially degenerate} if it is $h$--tangentially degenerate for
some $h$.\end{defn}

Remark that, if $h>0$, to say that $Z$ is $h$--tangentially degenerate
is equivalent to say that
$T_{Z,z}\cap Z$ has an irreducible component of dimension $h$.

We notice the following results:

\begin{prop}\label{tandeg} Let $Z\subseteq \pp^r$ be an irreducible,
non--degenerate, projective variety of
dimension $n$. One has:\pas
(a) $Z$ is $n$--tangentially degenerate if and only if $r=n$ and
$Z=\pp^n$;\pas
(b) if $n\geq 2$, then $Z$ is $(n-1)$--tangentially degenerate if and
only if either $r=n+1$ or $Z$
is a scroll over a curve.\end{prop}

\begin{proof} Part (a) is obvious. Part (b) is classical, and it can be
easily deduced from \cite
{ein}, Theorem 2.1 and Theorem 3.2.\end{proof}

\begin{lem} \label {tanproj} Let $X\subset\pp^r$, be an
irreducible, projective variety, with $r>2n+1$. Let $x\in X$ be
a general point and suppose that the projection of $X$ from $x$ is a
variety $X'\subset \pp^{r-1}$
which is not tangentially degenerate.  Then a general tangential
projection of $X$ is a
birational map to its image.\end{lem}

\begin{proof} If the conclusion
does not hold, then for $x,y\in X$ general points, the $(n+1)$--dimensional space
spanned by
$x$ and $T_{X,y}$ meets $X$ at another point $x'$. But then, by projecting $X$
from $x$, we get a variety
$X'\subset \pp^{r-1}$ such that the tangent space to $X'$ at a general
point $z$ meets $X'$ in
some other point $z'$. \end{proof}

\subsection{}\label {Scorzas} By Terracini's lemma and elementary considerations (see Proposition \ref {tandeg}), there are no defective curves. The study of defective surfaces goes back to
Palatini \cite {Pal1}, Scorza \cite {Scorza2}, Terracini \cite {Terr2},
whose classification result is the first complete one. In modern times, we
mention Dale \cite {Dale2} and Catalano--Johnson \cite {CJ}, whereas in
\cite {WDV} one finds the full classification of $k$--weakly defective
surfaces for any $k$, which includes Terracini's classification of 
defective surfaces. \par

The classification of $1$--defective threefoolds was taken up by Scorza in \cite
{Scorza}. The case of smooth $1$--defective threefolds
was also examined by Fujita \cite {Fuj} and Fujita--Roberts  \cite
{FujRob}. Scorza's classification has been revisited in 
\cite {ChCi}. We recall here the result:

\begin{thm} \label {scorza} Let $X\subset\pp^r$ be an irreducible,
non--degenerate, projective, $1$--defective threefold. Then $r\geq 6$ and
$X$ is of one of the following
types:

\begin {itemize}

\item [(1)] $X$ is a cone over a surface $S$;\pas
\item [(2)] $X$ sits in a $4$--dimensional cone over a curve $C$;\pas
\item [(3)] $r=7$ and $X$ sits in a $4$--dimensional cone over the Veronese
surface in $\pp^5$;\pas
\item [(4)] $X$ is either the $2$--Veronese embedding of $\pp^3$ in $\pp^9$
or a projection
of it in $\pp^r$, $r=7,8$;\pas
\item [(5)] $r=7$ and $X$ is a hyperplane section of the Segre embedding of
$\pp^2\times\pp^2$ in $\pp^8$.
\end{itemize}

Conversely, all the irreducible threefolds in the above list are
$1$-defective. \par 

Furthermore $\delta_1\leq 2$ and $m_1\leq 2$ and actually
$m_1=\delta_1=1$ unless we are in case (1) and
$X$ is either a cone over a curve or a cone over a Veronese surface in
$\p^5$, in which case $m_1=\delta_1=2$.
\end{thm}

\begin{rem} \label {remscorza}\rm It is useful to remark that, for a
$1$--defective threefold $X$ as in the various cases listed in Theorem \ref
{scorza} above, the general tangential $1$--contact locus is:

\begin {itemize}

\item [(1)] a reducible curve consisting of two general rulings of $X$,
unless $S$ is a $1$--weakly defective surface;
\item [(2)] a reducible surface consisting of two general fibres of the
projection of $X$ to $C$;
\item [(3)] an irreducible surface, the pull--back on $X$ of a general conic
of  the Veronese surface in $\pp^5$;\pas
\item [(4)] an irreducible conic, the image of a general line in
$\pp^3$;\pas
\item [(5)] an irreducible conic (\cite {ChCi}, Example 2.5).
\end{itemize}

Notice that:

\begin{itemize}

\item if the general tangential $1$--contact locus is a reducible
curve, then we are in case (1) and $X$, being a cone over a non--defective
surface, lies in $\p^r$, with $r\geq 7$;

\item if the general tangential $1$--contact locus is a reducible
surface, then we are either in case (2) or in case (1) and $X$ is the cone
over a $1$--weakly defective surface with reducible general tangential
$1$--contact locus. By looking at the classification of weakly defective
surfaces (see \cite {WDV}, Theorem 1.3), we see that also in this latter
case $X$ sits in a cone of dimension (4) over a curve and that $X$ sits
in $\p^r$, with $r\geq 7$;

\item if the general tangential $1$--contact locus is an irreducible
curve, then we are in case (4) or (5). Let $\mathcal G$ the family of 
tangential $1$--contact loci. This is a $4$-- dimensional family of conics. 
If $x\in X$ is a general point, we have therefore a $2$--dimensional family
$\mathcal G_x$ of conics in $\mathcal G$ passing through $x$. It is important to
notice, for future purposes, that the tangent lines to the conics in $\mathcal G_x$ 
fill up the whole  tangent space $T_{X,x}$. This is trivial in case (4), and
in case (5) it immediately follows by the discussion in \cite {ChCi}, Example 2.5. 

\end{itemize}
\end{rem}

In connection with the privious classification results, it is worth pointing out
the following:

\begin{prop}\label {rnc} Let $X\subset \p^ r$, $r\geq (k+1)n+k$, be an irreducible,
non--degenerate projective variety. Suppose that $X$ enjoyes the property that
$k+1$ general points of $X$ lie on a curve $\Gamma\subset X$ and such that $\dim(<C>)\leq 2k$. Then $X$ is $k$--defective. 
\end{prop}

\begin{proof}  Notice that $\Gamma$ is
not defective and therefore $S^{k}(\Gamma)=<\Gamma>$. Furthemore there are
infinitely many $(k+1)$--secant $\pp^ k$'s to $\Gamma$ and therefore to $X$,
containing the general point of $<\Gamma>$. This shows that there are infinitely
many $(k+1)$--secant $\pp^k$'s to $X$ containing the general point of $S^{k}(X)$,
thus $X$ is $k$--defective. 
\end{proof} 

\section {BASIC PROPERTIES OF DEFECTIVE THREEFOLDS} \label {defects}

From now on we will concentrate on defective threefolds. 
In this section we make a remark which, though not really needed in the
sequel, is, in our opinion, of some interest. 

We start by noticing that the basic tool for the proof of 
Theorem \ref {scorza} is the study of the first general
tangential projection $\tau_{X,1}: X\map X_1$. One studies separately the
cases in which $X_1$ is either
a curve or a weakly--defective surface. In these situations, one has
interesting information
about the general hyperplane section of $X$ which are important steps
toward the
classification. When $k>1$, the corresponding analysis is not
conclusive. However, one has
the following results:

\begin{prop} Let $X\subset \pp^r$ be an irreducible, non--degenerate,
projective, minimally $k$--defective
threefold. Assume that $m_k=2$, i.e. $X_k$ is a curve. Let $S$
be a general hyperplane
section of $X$. Then the general $k$--tangential projection $S_k$ of $S$
sits in a cone with a vertex
of dimension at most $k-2$ over the curve $X_k$.\par

Moreover:

\begin{itemize}

\item [(i)] when $r\ge 6k$, then $S$ is $(2k-1)$--defective;\pas
\item [(ii)] when $r\ge 6k-2$, then $S$ is $(2k-1)$--weakly defective.
\end{itemize}
\end{prop}

\begin {proof} Since $X$ is $k$--defective, we have that $S^k(X)$ is a
proper subvariety of $\p^r$, hence $k<r$. Let $p_0,...,p_k$ be general
points on $X$. Then $<p_0,...,p_k>$ is a proper subspace of $\p^r$. Let $H$
be a general hyperplane containing $<p_0,...,p_k>$. By
varying the points $p_0,...,p_k$ on $X$, then $H$ also varies and it turns
out to be a general hyperplane in $\p^r$. Let $S$ be the surface cut out by
$H$ on $X$. Then:

\begin{equation}\label{for1}
\dim(H\cap T_{X,P_1,\dots ,P_k})-\dim(T_{S,P_1,\dots ,P_k})\leq
k-1.\end{equation}

Notice that $S_k$, which is the projection of $S$ from
$T_{S,P_1,\dots ,P_k}$, is contained in a cone
$W$ over the projection $S'$ of $S$ from $H\cap T_{X,P_1,\dots ,P_k}$. We
may assume that
$W$ is of minimal dimension with the property of containing $S_k$. The
vertex $\Pi$ of $W$
is contained in the projection of $H\cap T_{X,P_1,\dots ,P_k}$ from
$T_{S,P_1,\dots ,P_k}$, hence,
by (\ref {for1}), we have $\dim(\Pi)\leq k-2$. Remark that $S'$ is also the
projection of $S$ from $T_{X,p_1,...,p_k}$. Hence $S'$ is contained in
$X_k$. Since $X_k$ is a curve, we have $S'=X_k$, proving the first part of
the assertion.\par

The tangent space to $S_k$ at a point $q$ intersects in a point the vertex
$\Pi$ of the cone $W$, because the projection of $S_k$ from
$\Pi$ is a curve (see Proposition \ref {varincon}). Let $q_1,\dots ,q_{k-1}$
be general points of
$S_k$. Then, by the minimality assumption on $W$,
the space $T_{S_k,Q_1,\dots ,Q_{k-1}}$ contains $\Pi$. Thus the general
$(k-1)$--tangential projection of $S_k$, which is the general
$(2k-1)$--tangential projection of $S$, certainly has dimension at most
$1$.  \par

Assume that $S$ is not $k$--defective, otherwise there is nothing to
prove. Then $S_k$ sits in
$\pp^s$, $s=r-3k-1$. From the previous argument, we get that the span of
$k$ general
tangent planes to $S_k$ has dimension at most $3k-2$. Thus, if $s\geq
3k-1$, i.e. if
$r\geq 6k$, then $S_k$ is
$(k-1)$--defective, hence $S$ is $(2k-1)$--defective. This
proves part (i).\par

If the span of $k-1$ general tangent planes to $S_k$ has dimension at most
$3k-5$ and if $s\geq
3k-4$, i.e. if $r\geq 6k-3$, then $S_k$ is $(k-2)$--defective, hence $S$ is
$(2k-2)$--defective and therefore also $(2k-1)$--defective. \par

Assume that the span of $k-1$ general tangent planes to $S_k$ has
dimension $3k-4$. Then the projection of $S_k$ from $T_{S_k,q_1,\dots
,q_{k-2}}$ is a surface and therefore
$T_{S_k,q_1,\dots ,q_{k-2}}$ cannot contain the vertex $\Pi$ of the cone
$W$. Since, as we saw, all spaces $T_{S_k,q_i}$ intersect $\Pi$ in a point,
by the minimality assumption on $W$ we have that 
$T_{S_k,q_1,\dots ,q_{k-2}}\cap \Pi$ has codimension $1$ in $\Pi$. Now the
projection $S_{2k-2}$ of $S_k$ from $T_{S_k,q_1,\dots ,q_{k-2}}$ is a
surface which sits in the cone $Z$ over
$X_k$ with vertex the projection of $\Pi$ from $T_{S_k,Q_1,\dots
,Q_{k-2}}\cap \Pi$, which is a
point. Hence $S_{2k-2}$ coincides with the cone $Z$. Moreover $S_{2k-2}$
spans a
$\pp^{s'}$ where $s'=s-3(k-2)=r-3k-1-3(k-2)=r-6k+5$. If $s'\geq 3$, i.e. if
$r\geq 6k-2$,
then $S_{2k-2}$ is $1$-weakly defective. Then $S$ is $(2k-1)$-weakly
defective (see \cite
{WDV}, Proposition (3.6)). \end{proof}

Similar arguments lead to the following:

\begin{prop} Let $X\subset \pp^r$ be an irreducible, non--degenerate,
projective, minimally $k$--defective threefold.  Assume that
$m_k=1$ and that $X_k$ is a developable surface. Let $S$ be a general
hyperplane section of $X$. Then the general
$k$-tangential projection $S_k$ of $S$ sits in a cone of dimension $k+1$
over the developable surface $X_k$. When $r\ge6k+1$, then $S$ is
$(2k-1)$-weakly defective.\pas \end{prop}

\section{BASIC PROPERTIES OF THE CONTACT LOCI}\label {2?}

In this section we study, mostly in the case of defective threefoolds, some
basic properties of the contact loci $\Gamma$ and $\Sigma$
defined in \S \ref {notaz}.\par

We start with the following:

\begin{prop}\label{mon}  Let $X\subset \pp^r$ be an irreducible,
non--degenerate, projective, $k$--defective variety. For a general choice
of $p_0,\dots,p_k\in X$ and a general choice of the hyperplane $H$
containing 
$T_{X,p_0,...,p_k}$, the contact loci
$\Gamma=\Gamma_{p_0,\dots,p_k}$ and
$\Sigma=\Sigma(H)=\Sigma_{p_0,\dots,p_k}(H)$ are equidimensional and
smooth at each of the points $p_0,...,p_k$. Furthermore either they are
irreducible or they consist of
$k+1$ irreducible components, one through each of the points
$p_0,\dots,p_k$.\end{prop}

\begin{proof} We prove the assertion for $\Gamma$. The proof for $\Sigma$
is quite
similar.\par

First of all, let us  move slightly the points $p_i$'s on some
component of $\Gamma$ to a new set of points $\{q_0,\dots,q_k\}$. Then
$q_0,\dots,q_k$ are also general points on $X$. Furthermore
$T_{X,p_0,\dots,p_k}$ contains the tangent spaces to $X$ at the points
$q_i$'s, so for dimension reasons, it coincides with
$T_{X,q_0,\dots,q_k}$. Hence $T_{X,q_0,\dots,q_k}$ is also tangent to $X$
along $\Gamma=\Gamma_{p_0,\dots,p_k}$. This tells us that, by the
generality of the points $p_0,\dots,p_k$, $\Gamma$ is
smooth at $p_0,...,p_k$, and therefore there is only one
irreducible component of $\Gamma$ through each of the points $p_0,\dots
,p_k$. Since $[p_0,...,p_k]$ is a general point of ${\rm Sym}^{k+1}(X)$, 
there is the monodromy action of the full symmetric group $\mathcal
S_{k+1}$ on $p_0,...,p_k$ as recalled in \S \ref {monodromy}. Hence we can
permute the points
$p_i$ as we like, and this implies that all components of
$\Gamma$ have the same dimension.\par

Assume now that there is a component of $\Gamma$ which contains more than
one of the points $p_i$'s, say $p_0$ and $p_1$. Again by monodromy, we can
let $p_0$ stay fixed and we can move $p_1$ to any one of the points $p_i$,
$i>1$. Then we see that also $p_0$ and $p_i$, $i>1$, stay on an irreducible
component of $\Gamma$. Since $p_0$ sits on only one irreducible component
of $\Gamma$, then this component has to contain all the points $p_0,\dots
,p_k$ and therefore it has to coincide with $\Gamma$. This proves the
proposition.\end{proof}

Since $\Gamma\subseteq \Sigma$, the next corollary immediately follows:

\begin {cor} Let $X\subset \pp^r$ be an irreducible, non--degenerate,
projective, $k$--defective variety of dimension $n$. If
$\epsilon_k=\gamma_k$, e.g. if either
$\epsilon_k=1$ or $\gamma_k=n-1$, then
the general contact locus $\Sigma$ coincides with the general tangential
contact locus $\Gamma$. \end{cor}

The next proposition, which will be useful later, tells us how tangential
contact loci behave under a general tangential projection:

\begin{prop}\label{reduccontact} Let $X\subset \pp^r$ be an irreducible,
non--degenerate,
projective, $k$--defective, but not $1$--defective variety of dimension $n$,
hence
$k\geq 2$.  Let
$X_1$ be the general
tangential projection of $X$. Then the general tangential $k$--contact
locus $\Gamma$ of $X$
and  the general tangential $(k-1)$--contact locus $\Gamma_1$ of $X_1$
have the same
dimension and $\Gamma$ is reducible if and only if $\Gamma_1$ is.
\end{prop}

\begin{proof} Since $X$ is not $1$--defective, we have $n_1=n$.
Hence, if $p_0\in X$ is a general point, the  projection $\tau:=\tau_{X,1}:
X\map X_1$ from
$T_{X,p_0}$ is generically finite.
If $p_1,\dots, p_k$ are general points of $X$, so that $\tau(p_1),\dots,
\tau(p_k)$ are
general points of $X_1$, then the image of the span
$T_{X,p_0,p_1,\dots,p_k}$ via $\tau$ is the
span $T_{X_1,\tau(p_1),\dots,\tau(p_k)}$. Hence $\Gamma$ maps onto
$\Gamma_1$ via $\tau$. The
generic finiteness of $\tau$ implies $\dim(\Gamma)=\dim(\Gamma_1)$.
Moreover if $\Gamma_1$ is
reducible then also $\Gamma$ is.\par

Conversely, suppose $\Gamma=\Gamma_{p_0,\dots ,p_k}$ is reducible. Then
$\Gamma=\cup_{i=0}^k\Gamma^{(i)}$,
where $\Gamma^{(i)}$ is the irreducile component containing the point
$i_i$. Set
$\Gamma^{(i)}_1=\tau(\Gamma^{(i)})$, so that $\Gamma^{(i)}_1$ is the
component of $\Gamma_1$
through the point $\tau(p_i)$. Since $p_1,\dots ,p_k$ are general points of
$X$, then
$\tau(p_1),\dots, \tau(p_k)$ are
general points of $X_1$. Therefore, once $p_1$, and thus
$\Gamma^{(1)}_1$, has been fixed, we can choose $p_2$ so that
$\tau(p_2)\notin \Gamma^{(1)}_1$, and so $\Gamma^{(1)}_1\not=
\Gamma^{(2)}_1$. Proceeding in this way we see that for a general
choice of  $p_1,\dots ,p_k$, the varieties
$\Gamma^{(1)}_1,
\dots , \Gamma^{(k)}_1$ are all distinct
varieties, proving that $\Gamma_1$ is reducible.\end{proof}

Now we restrict to the threefold case. One can easily detect when
$\Gamma$, and therefore
$\Sigma$, is a divisor, from the behaviour
of the tangential projection of $X$:

\begin{prop}\label{sigmasup} Let $X\subset \pp^r$ be an irreducible,
non--degenerate, projective,
minimally $k$--defective threefold. Then $\gamma_k=2$ if and only if either
one of the following holds:
\begin{itemize}

\item [(i)] $X_k$ is a curve;\pas
\item [(ii)] $X_k$ is a developable surface.
\end{itemize} \end{prop}

\begin{proof} We know
that $0<n_k<3$. Suppose that $\gamma_k=2$. Choose general points
$p_0,\dots ,p_k\in X$.
Let $\Gamma_0$ be the irreducible component of $\Gamma_{p_0,\dots ,p_k}$
containing $p_0$. Consider the
projection $\tau_{X,k}$ from the space $T_{X, p_1,\dots,p_k}$. Then
$\tau_{X,k}(\Gamma_0)$ is contained in $\tau_{X,k}(T_{X, p_0,\dots,p_k})$,
which is a projective space of dimension $n_k$. This shows that 
$\tau_{X,k}(\Gamma_0)$ is a proper subvariety of $X_k$, otherwise $X_k$
would be equal to the linear space $\tau_{X,k}(T_{X, p_0,\dots,p_k})$,
whereas we know that $X_k$ is a proper, non--degenerate subvariety of
$\p^{r_k}$ (see Proposition \ref {difetti}). Hence
$\tau_{X,k}(\Gamma_0)$ is either
a point or a curve. In the former case the fiber of $\tau_{X,k}$ at
$p_0$, which is a
general point of $X$, contains a divisor, hence $n_k=1$. In the
latter case, the
tangent plane to $X_k$ at its general point $\tau_{X,k}(p_0)$ meets $X_k$
along the curve $\tau_{X,k}(\Gamma_0)$, which passes through
$\tau_{X,k}(p_0)$. Hence $X_k$ is a
developable surface.\par

Conversely if either $n_k=1$ or $X_k$ is a developable surface, then
$\dim(\Gamma)=2$, because $\Gamma$ contains the pull--back, via
$\tau_{X,k}$, of the contact locus of
$X_k$ with its general tangent space. \end{proof}

Next, we look at the families $\mathcal G$ and $\mathcal S$ respectively
describing the
general contact loci $\Gamma=\Gamma_{p_0,\dots ,p_k}$ when $p_0,\dots
,p_k$ vary and by
$\Sigma=\Sigma(H)_{p_0,\dots ,p_k}$ when $p_0,\dots ,p_k$ and $H$ vary.
Recalling the definition of involution from \S \ref {invo}, we have:

\begin{lem}\label{involuzione} In the above setting, if $\gamma_k=2$
[resp. $\epsilon_k=2$], then the family of divisors
$\mathcal G$ [resp.
$\mathcal S$] is an involution of dimension $k+1$. Hence, if the general
member of $\mathcal G$ [resp. of
$\mathcal S$] is irreducible, then $\mathcal G$ [resp.
$\mathcal S$] is a linear system.\end{lem}

\begin{proof} To prove the first part of the assertion one argues as
in \cite {WDV}, pp. 172--173: though one refers to the surface case there,
the argument applies, word by word, to our situation. The final part
of the assertion follows by Theorem  \ref {CasHum}.\end{proof}

The following proposition gives more precise information about the
situation described in
Proposition \ref {sigmasup}:

\begin{prop}\label{minr} Let $X\subset \pp^r$ be an irreducible,
non--degenerate, projective,
minimally $k$--defective threefold. Assume that $\gamma_k=2$. Then:

\begin{itemize}

\item [(i)] when $X_k$ is a curve and $r>4k+2$, then the involution
$\mathcal G$ is composed with a
pencil;\pas
\item [(ii)] when $X_k$ is a developable surface, the same conclusion holds
provided $r>4k+3$.\end{itemize}
\end{prop}

\begin{proof} Assume that $r>4k+2$ and $X_k$ is a curve. By Proposition
\ref {difetti}, we know that $X_{k-1}$ is a threefold in $\pp^{r_{k-1}}$,
whose general tangential projection is
a curve. Notice that $r_{k-1}=r+4-4k>6$. Then by Theorem \ref {scorza} we
have that $X_{k-1}$ is a cone over a curve
with vertex a line and therefore two general tangent spaces to
$X_{k-1}$ have contact with $X_{k-1}$ along two planes. By pulling back to
$X$, we see that
$\Gamma_{p_0,\dots,p_k}$ has at least two irreducible components. It
follows from Theorem \ref{CasHum} that $\mathcal G$ is composed with a
pencil. This finishes case (i).\par

Similarly, assume that $X_k$ is a develobable surface and $r>4k+3$. By
Lemma 3.6 of
\cite{ChCi} a general hyperplane section $S'$ of $X_{k-1}$ is $1$--weakly
defective but not $1$--defective. By the classification of weakly
defective surfaces
(see \cite{WDV}, Theorem (1.3)), one gets that either:\pas

\begin{itemize}

\item [(a)] $S'$ is developable or\pas

\item [(b)] $S'$ is contained in a cone over a curve, with vertex along a
line or\pas

\item [(c)] $S'$ sits in $\pp^6$.\end{itemize}

Case (c) is excluded because $S'$ is non--degenerate in
$\pp^{r_{k-1}-1}$ and $r_{k-1}-1=r+3-4k>6$. Then only cases (a) and (b)
may occur and,
in both of them the tangential contact variety $\Gamma'$ of $S'$ with
the span of two
general tangent planes is reducible. \par

Now we claim that the general tangential contact locus of $X_{k-1}$ at
two general points is
reducible. This, as in case (i), leads to the assertion. As
for the claim, let
$q_1,q_2$ be two general points on $X_{k-1}$. Then $T_{X_{k-1},q_1,q_2}$
is a $\pp^6$,
because $\delta_k(X)=\delta_1(X_{k-1})=1$. Taking a general hyperplane
$H'$ through
$q_1,q_2$ and its section $S'$ with $X_{k-1}$, we have that
$T_{S',q_1,q_2}$ is a $\pp^5$,
because $S'$ is not $1$--defective. Hence $T_{S',q_1,q_2}=H'\cap
T_{X_{k-1},q_1,q_2}$. Let
$\overline \Gamma$ be the contact variety of $T_{X_{k-1},q_1,q_2}$ with
$X_{k-1}$. Of course the tangential contact locus
$\Gamma'$ for $S'$ is the intersection of $\overline \Gamma$ with $H'$.
Since
$\Gamma'$ is reducible
and by the genericity of $H'$, we see that $\overline \Gamma$ is
reducible. This proves the
claim and therefore the proposition.\end{proof}

\section{THE IRREDUCIBLE DIVISORIAL CASE} \label {irrediv}

Let $X\subset \pp^r$ be an irreducible, non--degenerate, projective,
minimally $k$--defective threefold with $k\geq 2$. In this section we
examine the case in
which $\gamma_k(X)=2$ and the general tangential $k$--contact
locus is irreducible. \par

By Proposition \ref {sigmasup} we know that either $X_k$ is a curve or
$X_k$ is a developable surface. Furthermore Proposition  \ref{minr} and
Proposition \ref {difetti} imply that $4k+2\leq r\le 4k+3$ and if $r=4k+2$
then $X_k$ is a curve.\par

Recall that, by Lemma
\ref {involuzione}, the family of divisors
$\mathcal G$ of tangential
$k$--contact loci is a linear system of dimension $k+1$ which is not
composed with a pencil.\par

Observe that, by construction, the linear system
$2\mathcal G$ (see \S \ref {linsyst}) is contained in the hyperplane linear
system
$\mathcal H$ of
$X$. In particular we have the linear equivalence relation $H \equiv
2\Gamma + \Delta$, with $\Delta$ effective.\par

We let:

$$\phi_\mathcal G: X\map Y\subset \p^{k+1}$$

\noindent be the rational map defined by $\mathcal G$. Since $\mathcal
G$ is not composed with a pencil, then $n:=\dim(Y)>1$. We set $d:=\deg(Y)$.
The classification is based on the classification of $Y\subset
\pp^{k+1}$ according to its dimension and degree. A straightforward
application of Theorem
\ref {Hilb} gives the following useful information:

\begin{lem} In the above setting, one has:

\begin{equation} \label{Castelnuovo} 4k+4\ge r+1\ge h_Y(2)\ge
\iota+k(n+1)-{\frac{d(d-3)} 2}+2\end{equation}
where $\iota=\iota(Y)=\min\{d-k+n-2,k+1-n\}$. Moreover if $d<2(k-n)+3$, the
equality holds in \eqref{Castelnuovo} if and ony if the general
curve section of $Y$ is projectively normal.\end{lem}

We will discuss separately the two cases $n=2,3$.

\begin{thm}\label {dimY3} In the above setting, assume that $n=3$. Then
$\iota\leq 3$,
$\Delta=0$ and we are in one of the following cases:
\begin{itemize}

\item [(1)] $\iota=0$ and $Y$ is a threefold of minimal degree in
$\pp^{k+1}$; then $X$ is contained in a cone over the $2$--uple embedding of
$Y$, with vertex either a point if $r=4k+2$ or a line if $r=4k+3$;\pas

\item [(2)] $\iota=1$, $k=3$, and $Y$ is a hypersurface of
degree $d\geq 3$ in $\pp^4$; then
either $r=4k+2=14$ and $X$ is the $2$-uple embedding of $Y$ or
$r=4k+3=15$ and $X$ is
contained in the cone with vertex a point over the $2$-uple embedding of
$Y$; \pas

\item [(3)] $\iota=1$ and $Y$ is a threefold of degree $k$ in $\pp
^{k+1}$ with curve sections of arithmetic genus $1$; then either
$r=4k+2$ and $X$ is the
$2$-uple embedding of $Y$ or $r=4k+3$ and $X$ is contained in the cone
with vertex a point
over the $2$-uple embedding of $Y$;\pas

\item [(4)] $\iota=1$ and $Y$ is a threefold of degree $k$ in $\pp ^{k+1}$
with curve sections
of genus $0$ which is either a cone with vertex a line over a smooth
rational curve of degree
$k$ in $\pp^{k-1}$ or it has a double line; then $r=4k+3$ and $X$ is the
$2$-uple embedding
of $Y$;\pas  

\item [(5)] $\iota=2$, $k=4$ and $Y$ is a threefold of degree $d\geq
5$ in $\pp ^5$,
contained in a quadric; then $r=4k+3=19$ and $X$ is the $2$-uple
embedding of $Y$;\pas

\item [(6)] $\iota=2$, $k\geq 4$ and $Y$ is a threefold of degree $k+1$ in
$\pp ^{k+1}$ with curve sections of
arithmetic genus 2; then $r=4k+3$ and $X$ is the $2$-uple embedding of
$Y$.\end{itemize}\end{thm}

\begin{proof} In the present case equation (\ref {Castelnuovo}) reads:

$$4k+4\ge r+1\ge h_Y(2)\ge \iota+4k+2$$
hence $0\leq \iota\leq 2$. Then $2\mathcal G$ is a subsystem of
$\mathcal H$ of codimension $\le 2$. It follows that $\Delta$ imposes at
most two
conditions to the hyperplanes of $\pp^{k+1}$. Since $\Delta$ is a
divisor, the only
possibility is that $\Delta=0$.\par

If $\iota=0$, then Theorem \ref {Hilb} yields $d=k-1$, i.e. $Y$ is
a threefold of minimal degree in $\pp^{k+1}$.  Since $2\mathcal G$ is a
subsystem of
codimension at most $2$ in $\mathcal H$, we are in case (1).\par

If $\iota=1$ then either $d=k$
or $k=3$ and $d\geq 3$. In the latter case we are in case (2). In
the former
case if $h_Y(2)=4k+3$, then the curve section of $Y$
is linearly normal of degree $k$ and
arithmetic genus $1$ (see Theorem \ref {Hilb}) and we are in
case (3). If
$h_Y(2)=4k+4$ then, again by Theorem \ref {Hilb}, the curve section of $Y$
cannot be linearly normal of degree $k$ and
arithmetic genus $1$, hence it has arithmetic genus $0$ and it is
therefore smooth. By Lemma
\ref {RNS} we are in case (4). \par

If $\iota=2$, then all the inequalities in (\ref {Castelnuovo}) are
equalities, hence $X$ is the
$2$--uple embedding of $Y$. We have either $k=4$ and $d\geq 5$ or
$d=k+1$ and $k\geq 4$. If
$k=4$, the equalities in (\ref {Castelnuovo}) imply that $Y$ is
contained in a quadric of
$\pp^5$ and we are in case (5). If $d=k+1$, then $Y$ has degree
$k+1$. Since
equalities hold in (\ref {Castelnuovo}), by Theorem \ref {Hilb},
a general
curve section of $Y$ is linearly normal of degree $k+1$, with arithmetic
genus $2$ and we are in case (6). \end{proof}

\begin{ex}\label{exx} We show that if $X$ is as in any one of the
previous cases, then it is minimally
$k$--defective. We may assume $k\geq 2$. \medskip

\rm \noindent (1) Assume that $X\subset \p^r$ is contained in a cone $W$
over the $2$-uple embedding
$Y'\subset \pp^{4k+1}$ of a threefold $Y$ of minimal degree in
$\pp^{k+1}$ with vertex $L$ a
line or a point. Accordingly one has $4k+3\geq r \geq 4k+2$. Remark that by definition
$X$ projects onto $Y'$ from $L$. \par 

In Example \ref {ex'}
below it is proved that $Y'$ is
$k$--defective and $s^{(k)}(Y')=4k$. By Proposition \ref {proiezione},
part (i), one has $s^{(k)}(X)\leq
s^{(k)}(Y')+\dim(L)+1= 4k+\dim(L)+1<\sigma^{(k)}(X)$, thus
$X$ is $k$--defective. \par

Let us prove that $X$ is minimally $k$--defective. In Example \ref {ex'}
below we will see that $Y'$ is minimally $(k-1)$--defective and that
$s^{(k-1)}(Y')=4k-2$. Moreover the general $(k-1)$--contact locus of $Y'$
is an irreducible curve. Suppose, by contradiction, that $X$ is
$(k-1)$--defective. Then by Proposition \ref {proiezione}, part
(i), we have $4k-2=s^{(k-1)}(Y')\leq s^{(k-1)}(X)<\sigma^{(k-1)}(X)=4k-1$,
thus $s^{(k-1)}(X)=4k-2$. By Proposition \ref {proiezione}, part
(iii), $X$ it would be, as well as $Y'$, minimally
$(k-1)$--defective. Moreover the general
$(k-1)$--contact locus of $X$ projects from $L$ to the general
$(k-1)$--contact locus of $Y'$, hence it is an irreducible curve. However
we will see at the beginning of \S  \ref{irredcurv} below that the maximum for the embedding
dimension $r$ of a minimally $(k-1)$--defective threefolds whose general 
$(k-1)$--contact locus is an irreducible curve is $r\leq 4k+1$. This
gives a contradiction, which proves that $X$ is not $(k-1)$--defective.\par

If $L$ is a point, then $r=4k+2$ and by Proposition \ref {difetti},
part (v), we have $\delta_k(X)=1$, $X_k$ is a plane curve and
$s^{(k)}(X)=4k+1$. \par

If $L$ is a line, then $r=4k+3$.  By Proposition \ref {minr}, part (i), we
have $n_k=2$. Then Proposition
\ref {difetti}, part (iii), yields $\delta_k(X)=1$ hence $s^{(k)}(X)=4k+2$.
Actually we see that in this case $X_k$ is a cone in $\p^3$. Indeed, by
projecting from a general point of $L$ one has to go back to the
previous situation.\medskip

\noindent (2) Let $Y$ be a threefold in $\pp^4$ of degree $d\geq
3$. Let $X:=Y'$ be the $2$--uple embedding of $Y$ in $\p^{14}$. If $p$ is a
point of $Y$ we abuse notation and denote by $p$ also the corresponding
point on $Y'$.

 We have
$\sigma^{(3)}(X)=14$. Let $p_0,...p_3$ general points on $Y$.
Since there is a quadric in $\pp^4$
singular at $p_0,...p_3$, namely the hyperplane $<p_0,...p_3>$ counted
twice, then
$T_{Y',p_0,...,p_3}$ is contained in a hyperplane. This proves that $Y'$ is
$3$--defective. \par

Let us prove that $Y'$ is minimally $3$--defective. Certainly $Y'$ is not
$1$--defective. Indeed let $p_0,p_1$ be general points on $Y$. There are
hyperplane sections of $Y'$ which have isolated
singularities at $p_0,p_1$, namely the
hyperplane sections of $Y$ corresponding to the quadrics in $\p^4$
singular along the line $<p_0,p_1>$. Suppose that $Y'$ is
$2$--defective and therefore minimally $2$--defective. Let $p_0,p_1,p_2$ be
general points on
$Y$. There are hyperplane sections of $Y'$ tangent at $p_0,p_1,p_2$
and having an irreducible singular locus of dimension $1$, namely the
hyperplane sections of $Y$ corresponding to the quadrics in $\p^4$
singular along the plane $<p_0,p_1,p_2>$. Then, as we will seee at the beginning of
\S  \ref{irredcurv} below, the embedding dimension of $Y'$ should be bounded above by $13$, a
contradiction.

By Proposition \ref {difetti}, part (v), we have $\delta_3(Y')=1$,
$s^{(3)}(Y')=13$, $n_3(Y')=1$ and $Y'_3$ is a plane curve.\par

Assume now that $X$ is contained in a cone with
vertex a point $v$ over $Y'$, so that $r=15$ and that $X$ maps to $Y'$ from
$v$. By Proposition \ref {proiezione}, part (i), we have $s^{(3)}(X)\leq
s^{(3)}(Y')+1=14$, hence $X$ is
$3$--defective. By the same Proposition \ref {proiezione}, part (ii), since
$Y'$ is minimally $3$--defective, also $X$ is minimally $3$--defective.\par

By Proposition \ref {minr}, part (i), we
have $n_k=2$. Then Proposition
\ref {difetti}, part (iii), yields $\delta_k(X)=1$ hence
$s^{(k)}(X)=14$.\medskip

\noindent (3) Assume that $X$ is contained in a cone with vertex a point (or
the empty set) over the
$2$-uple embedding $Y'$ of a threefold $Y$ of degree $k$ in $\pp^{k+1}$
with linearly
normal curve section of arithmetic genus 1. Then $Y'$ sits in
$\pp^{4k+2}$ (see Theorem \ref {Hilb}). In \cite{WDV}, Example 4.7, it is
proved that $Y'$ is $k$--defective and $s^{(k)}(Y')\leq 4k+1$. An
argument completely similar to the one of the previous Example (2) shows
that $Y'$ is minimally $k$--defective. \par

If $X=Y'$ then Proposition \ref {difetti}, part (v), yields
$\delta_k(X)=1$, $n_k(X)=1$, $s^{(k)}(Y')=4k+1$.\par

Assume $X\not= Y'$, hence $X\subset \p^{4k+3}$. Then Proposition \ref
{proiezione}, part (i), yields $s^{(k)}(X)\leq s^{(k)}(Y')+1=4k+2$, thus
$X$ is $k$--defective. Again Proposition \ref {proiezione} implies it is
minimally $k$--defective. Then by Proposition \ref {minr}, part (i), we
have $n_k=2$. Finally Proposition
\ref {difetti}, part (iii), yields $\delta_k(X)=1$, hence
$s^{(k)}(X)=4k+2$.\medskip

\noindent (4) There are threefolds $Y$ of degree $k$ in $\pp^{k+1}$ with
smooth curve sections of
genus $0$ and a singular line. They are described in Proposition \ref
{ratsing}. By Lemma \ref {RNS}, the $2$--uple embedding $X$ of $Y$ sits in
$\p^{4k+3}$. By arguing as in Example 4.5 or 4.7
of \cite {WDV}, one sees that $X$ is $k$--defective. 
Then an argument completely similar to the one of the previous Example (2)
shows that $X$ is minimally $k$--defective. As usual we find
$\delta_k(X)=1$, $n_k=2$ and
$s^{(k)}(X)=4k+2$.\medskip

\noindent (5) Let $X\subset \p^{19}$ be the $2$-uple embedding of a
threefold $Y$ contained in a unique quadric
of $\pp^5$. If $p$ is a point of $Y$ we abuse notation and denote by $p$
also the corresponding point on $X$.\par

One has $\sigma^{(4)}(X)=19$. However
consider five general points $p_0,\dots,p_4$ of $Y$. Since there is a
quadric singular at any five points of $\pp^5$, i.e. the
double hyperplane $<p_0,...,p_4>$, then there is a hyperplane section of
$X$ tangent at $p_0,...,p_4$, hence $X$ is
$4$--defective. The usual kind of arguments show that $X$ is minimally
$4$--defective and that $\delta_4(X)=1$, $n_4=2$ and
$s^{(k)}(X)=18$. \medskip

\noindent (6) Let $X$ be the $2$--uple embedding of a threefold $Y$ of
degree $k+1$ in $\pp^{k+1}$ with curve sections of arithmetic genus $2$. 
Then $X$ sits in $\p^{4k+3}$ by Thoerem \ref {Hilb}. Moreover $X$ is
$k$--defective. This follows by arguing as in
\cite{WDV}, Example 4.7 in the usual way. Again
one has $\delta_k(X)=1$, $n_k=2$ and
$s^{(k)}(X)=4k+2$.\medskip

Notice that for all the above examples $X\subset \p^r$, one has
$s^{(k)}(X)=r-1$. Then, by Proposition \ref {crescente}, 
$s^{(k+1)}(X)=r$, i.e. $X$ is not $h$--defective for any $h\not=
k$.\medskip

It is worth adding a few words about the
existence of smooth threefolds of the above types. It is obvious that this can
happen for threefolds of types (1), (2) and (5). Threefolds of types (3) and (6)
can be smooth, but only for finitely many values of $k$ (this follows, for
instance, from the results in \cite {hart}; more specifically, for the genus $1$
case, see \cite
{Murre} for the genus $1$ case). Threefolds of type (4) can never be smooth. 

\end{ex}

Next let us turn to the case $n=\dim(Y)=2$. First we prove a lemma on
the dimension of the system $2\mathcal G$.

\begin{lem}\label{littlebound} In the above setting, if $n=2$ then
$\dim(2\mathcal G)\le r-2$. If the equality holds, then $\Delta=0$.\end{lem}

\begin{proof} Notice that $X$ is not a cone because it is $k$--mimimally
defective, with $k>1$.\par

 The
system $2\mathcal G$ maps $X$ to a surface $Y$, then it cannot coincide
with $\mathcal H$. Since
$\mathcal H$ has no fixed parts, then $\dim(2\mathcal G)<r$.\par

 If $\dim(2\mathcal
G)=r-1$, then, as in Proposition \ref{dimY3}, we see that $X$ is
contained in the cone, with
vertex a point, over $Y$, which is a
surface. Then $X$ would be the cone over $Y$, a contradiction. \par

If $\dim(2\mathcal G)= r-2$, then $X$ sits in a cone over $Y$ with
vertex a line $\ell$. The divisor
$\Delta$ has to be contained in $\ell\cap X$, hence $\Delta=0$. \end{proof}

Now we are ready for the proof of the classification theorem in the case
$n=2$:

\begin{thm}\label{dimY2} In the above setting, assume $n=2$. Then
$r=4k+3-i$, $i=0,1$ and one of the following cases occurs:\pas

\begin{itemize}

\item [(1)] $X$ is contained in a cone with vertex a space of dimension
$k-i$ over the $2$--uple embedding of a
surface $Y$ of minimal degree in $\pp^{k+1}$;\pas

\item [(2)] $k=2$, $i=0$ and $X$ sits in a cone with vertex a line over the
$2$--uple embedding of a surface $Y$ of
$\pp^3$ with $\deg(Y)\geq 3$;\pas

\item [(3)] $k\geq 3$, $i=0$ and $X$ sits in a cone with vertex of dimension
$k-1$ over the
$2$--uple embedding of a surface $Y$ of degree $k+1$ in $\pp^{k+1}$ with
curve sections of arithmetic
genus $1$.

\end{itemize}

All threefolds in the above list are minimally $k$--defective with
$\delta_k(X)=1$,
$s^{(k)}(X)=r-1=4k+3-i$, $n_k=2-i$, and are not $h$--defective for any
$h>k$. \end{thm}

\begin{proof} Recall that $r=4k+3-i$ with $i=0,1$. Formula (\ref
{Castelnuovo}) and Lemma \ref {littlebound} give:

$$4k+2-i\geq (r+1)-2=r-1\geq h_Y(2)\geq 3k+3+\iota$$

\noindent where $\iota=\min(d-k,k-1)$.  Then, by writing
$h_Y(2)=3k+3+q$, we see that $X$ is contained in a cone $W$ with vertex
of dimension $s=r-3k-q-3=k-q-i$
over the $2$-uple embedding $Y'$ of a non degenerate surface
$Y\subset\pp ^{k+1}$, such that
$h_Y(2)=3k+3+q$. By Proposition \ref {conidif}, $X$ is certainly
$(s+1)$--defective. By the minimality assumption we must have $k\leq
s+1\leq k-q-i+1$, thus $q+i\leq 1$. \par

If $q=0$, by Theorem \ref {Hilb},  
$Y$ is a surface of minimal degree in $\pp^{k+1}$. In this case $Y'$ is
minimally $k$--defective (see Theorem (1.3) of \cite {WDV}) and, 
Proposition \ref {conidif} implies that $X$ is minimally $k$--defective
too. One has $s^{(k)}(X)=s^{(k)}(Y')+k-i+1=4k+2-i$, hence $\delta_k(X)=1$
and $n_k(X)=2-i$.\par 

If $q=1$, then $i=0$. Taking into account Theorem \ref {Hilb}, we see that
only the following cases may occur:

\begin{itemize}

\item $k=2$ and $Y$ is a surface in $\pp^3$ of degree $d\geq 3$
and we are in case (2);

\item $k\geq 3$ and $Y$ is a surface of degree $k+1$ in
$\pp^{k+1}$ with curve sections of
arithmetic genus $1$, and we are in case (3).\end{itemize}

In both the above cases $Y$ is not $k$--defective (see Theorem (1.3) of
\cite {WDV}), and therefore by, Proposition \ref {conidif}, $X$ is
minimally $k$--defective. \end{proof}

\begin {rem}\label {liscezza} \rm (1) Observe that in the last statement, in cases
(1) and (3) we cannot specify whether $\Delta\equiv
H-2\Sigma$ is zero or not. When it is not, then $X$ meets the vertex of
the cone along a divisor. In
case (2) instead, $\Delta=0$ by Lemma \ref {littlebound}. \medskip

\noindent (2)  A few words about the existence of smooth threefolds  in the list
of Theorem \ref {dimY2}. For
example, consider case (1). Give an embedding whatsoever of a smooth $Y$ surface
of minimal degree in $\p^{k+1}$ in a $\p^{k-i}$.  This is certainly possible if
$k-i>4$. Call $Y'$ the surface we get in this way. Then join any point of
$Y'$ with the corresponding point of $Y$. The resulting ruled threefold is
smooth.\par

With the same idea one produces smooth threefolds of type (3) in Theorem \ref
{dimY2} as soon as $k\geq 6$. However note here that such a construction
works only for $k\leq 8$, since, as it is well known, there are
no smooth surfaces of degree $k+1$ in $\pp^{k+1}$ with curve sections of arithmetic
genus $1$ as soon as $k\geq 9$. \par
\end{rem} 

\section{THE REDUCIBLE CASE} \label {red}

Let $X\subset \pp^r$ be an irreducible, non--degenerate, projective,
minimally $k$--defective threefold with $k\geq 2$. In this section we
examine the case in
which the general element of the family $\mathcal G$ is reducible. By
Proposition \ref{mon},
we know that in this case $\Gamma=\Gamma_{p_0,\dots ,p_k}$ has exactly $k+1$
components, one passing through each of the points $p_i$,
$i=0,...,k$.\par\smallskip

The following proposition immediately follows by Proposition \ref
{reduccontact},  Remark
\ref {remscorza}:

\begin{prop} \label {1stclass} Let $X\subset \pp^r$, be an irreducible,
non--degenerate,
projective, minimally $k$--defective threefold with $k\geq 2$. Suppose
that the general
tangential $k$--contact locus of $X$ is reducible. One has:\pas

\begin{itemize}

\item [(1)] if $\gamma_k=1$, then $r\geq 4k+3$, $\delta_k=1$, $n_k=2$ and
the general
$(k-1)$--tangential projection $X_{k-1}$ of $X$ is
a cone over a surface which is not $1$--weakly defective;\pas

\item [(2)] if $\gamma_k=2$ and  $n_k=2$, then $r\geq 4k+3$, $\delta_k=1$
and the general
$(k-1)$--tangential
projection $X_{k-1}$ of $X$ is a threefold contained in a
$4$-dimensional cone over a
curve;

\item [(3)] if $\gamma_k=2$ and
$n_k=1$ then $r\geq 4k+2$ and the general $(k-1)$--tangential
projection
$X_{k-1}$ of $X$ is a cone
with vertex a line over a curve.
\end{itemize}

Conversely, threefolds like in case (1), (2) and (3)  above are minimally
$k$--defective with reducible general
tangential $k$--contact locus.

\end{prop}

We can now specify case (1) of Proposition \ref {1stclass}:

\begin{prop} \label{proj} Let $X\subset \pp^r$ be an irreducible,
non--degenerate projective threefold. Then the following are equivalent:\par

\begin {itemize}

\item [(i)] the general $(k-1)$--tangential projection $X_{k-1}$
of $X$ is a cone over a surface which is not $1$--weakly defective;

\item [(ii)] $X$ sits in a cone of dimension $k+2$, and not smaller, over a surface which is not $k$--weakly defective. \par

\end{itemize}

\end{prop}

\begin{proof} Supppose the general $(k-1)$--tangential projection $X_{k-1}$
of $X$ is a cone over a non--developable surface $S$, hence  $X$ is
minimally $k$--defective. Part (i) is trivially true for $k=1$. So we
assume $k\geq 2$ and proceed by induction on $k$.\par

Let $p_0\in X$ be a general point and let $\tau:=\tau_{X,1}: X\map X_1$ be
the tangential projection from $T_{X,p_0}$. By Proposition \ref {1stclass} 
and Proposition \ref {reduccontact}, $\gamma_{k-1}(X_1)=1$ 
and the general tangential $(k-1)$--contact locus of $X_1$ is a reducible curve.
\par

By induction, we know that $X_1$ is
contained in a cone over a surface $S_1$ with vertex $\Pi:=\Pi_{p_0}$ of dimension $k-2$.  Notice that the
general tangent space to $X_1$ meets $\Pi$ at one point and projects from
$\Pi$ to the general tangent space to $S_1$ (see Proposition
\ref {varincon}). \par

Call $\pi:X \map S_1$
the composition of $\tau$ with the projection of $X_1$ from $\Pi$. For a general point $p\in X$, call $Z_p$ the intersection of
$X_1$ with $<\Pi,\pi(p)>$, which is the general fibre of the projection of $X_1$ from $\Pi$. Hence $Z_p$ is a curve. \par

Take general points
$p_1,\dots ,p_k\in X$. Terracini's lemma and Proposition \ref {conidif}
imply that  $T_{X_1,\tau(p_1),\dots ,\tau(p_k)}=
<\Pi,T_{S,\pi(p_1),\dots ,\pi(p_k)}>$, hence $T_{X_1,\tau(p_1),\dots
,\tau(p_k)}$ is tangent to $X_1$ along the curves
$Z_{p_i}$, $i=1,\dots ,k$. Furthermore $T_{X_1,\tau(p_1),\dots
,\tau(p_k)}$ is tangent to $X_1$ along the pull--back of the general tangential $(k-1)$--contact locus to $S_1$. Since, as we saw, 
the general tangential $(k-1)$--contact locus of $X_1$ is a reducible curve, 
we have that $S_1$ is not $(k-1)$--weakly defective. \par

Notice that, if we move $p_1$ to a new general
point $p'_1$, then 
$T_{X_1,\tau(p'_1),\dots ,\tau(p_k)}=
<\Pi,T_{S,\pi(p'_1),\dots ,\pi(p_k)}>$ is again tangent to $X_1$ along the curves $Z_{p_i}$, $i=2,\dots ,k$. 

Call $C_i$ the pull--back to $X$ of $Z_{p_i}$ via $\tau$. One has that
$T_{X,p_0,p_1,\dots,p_k}$ is tangent to $X$ along a curve which contains
all the $C_i$'s. In
other words $\Gamma_{p_0,p_1,\dots,p_k}\supset C_1\cup\dots\cup C_k$. By
what we saw above, if we move the point $p_1$ to a new general point $p'_1$,
then $T_{p_0,p'_1,p_2\dots,p_k}$ is again tangent along the
curves $C_2,\dots ,C_k$. Therefore also by
moving $p_0$ to some other point $p_0'\in X$, again
$\Gamma_{p_0',p_1,\dots,p_k}$ contains
the curves $C_i$, $i=1,\dots ,k$.\par

By equation \eqref {terracini'}, we know that $h(\Gamma)\leq
2(1+k)-1=2k+1$, so $\dim(<\Gamma>)\leq 2k$. On the
other hand we claim that:

\begin{equation}\label{claim} \dim(<\Gamma>)\geq
k-1+\dim(<C_i>)\end{equation}
Otherwise there is an $h<k$ such that $\Lambda:=<C_1,\dots ,C_h>$ contains
all curves $C_i$,
with $h<i\leq k$. Then, by the genericity of the points $p_1,\dots,p_k$,
the space
$\Lambda$ would contain the whole of
$X$, but this is impossible,
since:

$$\dim(<\Lambda>)\leq\dim (<\Gamma>)\leq 2k<4k+2\leq r.$$

\noindent By (\ref {claim}) it follows that
$\dim(<C_i>)\leq k+1$. \par

Set $L:=L_{p_0}=<T_{X,p_0},\Pi>$. By Proposition  \ref {conidif}, $X$
cannot lie on a cone on a surface with vertex $\Pi$. Thus $\dim(L)=k+2$. We
have:

$$\dim (<Z_{p_i}>) = \dim(<C_i>)-\dim(<C_i>\cap T_{X,p_0})-1 $$
$$\dim (<Z_{p_i}>\cap \Pi) =
\dim( <C_i>\cap L)-\dim(<C_i>\cap T_{X,p_0})-1.$$
Since $\dim (<Z_{p_i}>\cap \Pi)=\dim (<Z_{p_i}>)-1$ because $Z_{p_i}$
projects from $\Pi$ to a
point of $S$, we have that also  $H_i:=<C_i>\cap L$ is a hyperplane in
$<C_i>$, $i=1,\dots ,k$.\par

Now move the point $p_0$ to a new point $p_0'\in X$, which gives us a
new vertex
$\Pi':=\Pi_{p'_0}$ and a new space $L':= L_{p'_0}=<T_{X,p'_0},\Pi'>$.
However, since
$\Gamma_{p_0',p_1,\dots,p_k}$  contains the curves $C_i$, $i=1,\dots ,k$,
we have that also
$<C_i>\cap L'$ is a hyperplane in $<C_i>$, $i=1,\dots ,k$.\par

Now $\dim(<L,L'>)\leq 2k+5< 4k+2\leq r$, so that $C_i$ cannot be
contained in
$<L,L'>$, for a general choice of $p_i$, otherwise $X$ would be contained
in $<L,L'>$. Hence the hyperplane
$H_i$ of
$<C_i>$ is fixed, as $p_0$ moves. Then the space $<H_1,\dots, H_k>$ does
not depend on
$p_0$, i.e. it is contained in $L_p$ for a general point $p\in X$. \par

Consider now
the linear space $H=\cap_{P\in U}L_P$, where $U$ is a suitable dense open
subset of $X$. Notice that
$L_p$ has to vary when $p$ varies in $X$,
because $T_{X,p}\subset L_p$. Hence $\dim(H)\leq \dim(L_p)-1=k+1$.
Moreover $H$ contains all
the hyperplanes $H_i$ of $<C_i>$ so it meets the tangent line of $C_i$
at $P_i$, which is a
general point of $X$. Thus $H$ intersects all the tangent spaces to
$X$. We conclude by Proposition \ref {varincon} that $X$ projects from $H$
either to a curve or to a surface. \par

If $\dim(H)\leq k-1$, then by Proposition \ref {conidif} and since $X$ is minimally $k$--defective, we have that 
$\dim(H)= k-1$ and $X$ projects from $H$ to a surface $S$ which is not $k$--weakly defective, proving (ii).

Assume 
$\dim(H)\geq k$. Then since for a general point $p\in X$ we have $H\cup
T_{X,p}\subset L_p$ and
$\dim(L_p)=k+2$ then $H$ would meet $T_{X,p}$. If $\dim(H)=k$ then
$\dim(H\cap T_{X,p})\geq 1$ and,
by Proposition \ref {varincon}, we would have that $X$ sits in a cone with vertex
$H$ over a curve. Then, by applying Proposition \ref {conidif}, we see that $X$ would be $(k-1)$--defective, a contradiction.  If $\dim(H)=k+1$ then $\dim(H\cap T_{X,P_0})\geq
2$ and, by Proposition \ref {varincon}, $X$ would be
contained in a subspace of dimension $\dim(H)+1=k+2$ a contradiction.\par

We have thus completed the proof of the fact that (i) implies (ii). \par

Let us assume now that (ii) holds. Then, by Proposition \ref {conidif}, $X$ is minimally $k$--defective and, since $S$ is not $k$--weakly defective, then $\gamma_k(X)=1$ and the general tangential $k$--contact locus is a reducible curve. Then,  by Proposition \ref  {1stclass}, (i) holds.
\end{proof}

The following proposition is proved with a similar argument, but
with a slight modification:

\begin{prop} \label{proj'} Let $X\subset \pp^r$ be
an irreducible, non--degenerate projective threefold. Let $k\geq 2$ be an integer. Then the
following are equivalent:\par

\begin {itemize}

\item [(i)] the general $(k-1)$--tangential projection $X_{k-1}$
of $X$ sits in a $4$--dimensional cone over a curve;

\item [(ii)]  $X$ sits in a cone of dimension $2k+2$, and not smaller, over a curve. \par

\end{itemize}\end{prop}

\begin{proof}. Let us prove that (i) implies (ii).  The assertion trivially
holds for $k=2$. So we assume
$k\geq 3$ and proceed by induction on $k$.\par

Let $p_0\in X$ be a general point and let $\tau:=\tau_1: X\map X_1$ be
the tangential
projection from $T_{X,p_0}$. By induction, we know that $X_1$ is
contained in a cone over a curve
$C$ with vertex $\Pi:=\Pi_{p_0}$ of dimension $2k-2$.  Call $\pi:X\map C$
the composition of $\tau$ with the
projection of $X_1$ from $\Pi$. For a general point $p\in X$, call $Z_p$ the intersection of
$X_1$ with the span
$<\Pi,\tau(p)>$, which is the fibre of the projection of $X_1$ from $\Pi$.
Hence $Z_Q$ is a surface.\par

Take general points
$p_1,\dots ,p_k\in X$. One
has:

$$\dim(T_{C,\pi(p_1),\dots ,\pi(p_k)})\leq 2k-1$$
hence:

$$\dim(<\Pi,T_{C,\pi(p_1),\dots ,\pi(p_k)}>)\leq 4k-2$$
and $T_{X_1,\tau(p_1),\dots ,\tau(p_k)}\subseteq
<\Pi,T_{C,\pi(p_1),\dots ,\pi(p_k)}>$, so that:

$$\dim(T_{X_1,\tau(P_1),\dots ,\tau(P_k)})\leq 4k-2.$$
Since $\delta_k(X)=1$ (see  {1stclass} ) the equality has to hold in the previous inequality, which implies
$T_{X_1,\tau(p_1),\dots ,\tau(p_k)}=<\Pi,T_{C,\pi(p_1),\dots ,\pi(p_k)}>$.
Furthermore
$T_{X_1,\tau(P_1),\dots ,\tau(P_k)}$ is tangent to $X_1$ along the surface
$Z_{P_i}$,
$i=1,\dots ,k$. \par

Call $S_i$ the pull--back to $X$ of $Z_{P_i}$ via $\tau$. It follows that
$T_{X,p_0,p_1,\dots,p_k}$ is tangent to $X$ along a surface which
contains all the $S_i$'s. In
other words $\Gamma_{p_0,p_1,\dots,p_k}\supset S_1,\dots,S_k$. In
particular the
tangential $k$--contact locus of $X$ is reducible. As in the proof of
Proposition \ref
{proj} we see that by
moving $p_0$ to some other point $p_0'\in X$, then
$\Gamma_{p_0',p_1,\dots,p_k}$ also contains
the surfaces $S_i$, $i=1,\dots ,k$.\par

By equation \eqref {terracini'}, one has $\dim(<\Gamma>)\leq 3k+1$. On the
other hand, as in the proof of proposition \ref {proj}, one has
$\dim(<\Gamma>)\geq
k-1+\dim(<S_i>)$. Hence
$\dim(<S_i>)\leq 2k+2$. \par

Set $L:=L_{P_0}=<T_{X,p_0},\Pi>$, which is a linear space of dimension
$2k+2$. Again one
sees that  $H_i:=<S_i>\cap L$ is a hyperplane in $<S_i>$, $i=1,\dots ,k$.\par

Now move the point $p_0$ to a new point $p_0'\in X$, which gives us a new vertex
$\Pi':=\Pi_{p'_0}$ and a new space $L':= L_{p'_0}=<T_{X,P_0'},\Pi'>$. However, since
$\Gamma_{p_0',p_1,\dots,p_k}$ contains $S_i$, $i=1,\dots ,k$, we have that also
$<S_i>\cap L'$ is a hyperplane in $<S_i>$, $i=1,\dots ,k$.\par

Now we claim that $\dim(<L,L'>)<r$. In fact, $<T_{X,p_0},T_{X,p'_0}>$ contains
$T_{X_1,\tau(p'_0)}$ and $\Pi$ intersects this space along a line. Hence
$\dim(<T_{X,p_0},T_{X,p'_0},\Pi>)\leq 2k+4$. For analogous reasons,  $\Pi'$
intersects  $<T_{X,p_0},T_{X,p'_0}>$  at least along a line. Hence 
$\dim(<L,L'>)\leq 4k+1<r$. \par

Since $<L,L'>$ is a proper subspace of $\p^ r$, the general surface $S_i$ cannot be
contained in it. Hence the hyperplane $H_i$ of $<S_i>$ is fixed, as $p_0$ moves.
Then the space $<H_1,\dots, H_k>$ does not depend on $p_0$, i.e. it is contained in
$L_p$ for all points $p\in X$. \par

Now one concludes exactly as in the proof of Proposition \ref {proj}.

Let us prove now that (ii) implies (i). Suppose $X$ lies in a cone with vertex
$\Pi$ of dimension $2k$ over a curve $C$. By Proposition
\ref {conidif'},
$X$ is minimally $k$--defective, hence $X_{k-1}$ is a threefold. Since the general 
tangent space to $X$ meets $\Pi$ along a line, in the $(k-1)$--tangential $\Pi$
projects to a line $L$, and $X_{k-1}$ sits in the
cone with vertex $L$ over $C_{k-1}$.

\end{proof}

With similar arguments, using Proposition
\ref {conidif'}, one proves the following result:

\begin{prop} \label{proj''} Let $X\subset \pp^r$ be an irreducible, non--degenerate
projective threefold. Then the following are equivalent:\par

\begin {itemize}

\item [(i)] the general $(k-1)$--tangential projection $X_{k-1}$ of $X$ is a cone
over a curve;

\item [(ii)] $X$ sits in a cone of dimension $k+2$ over a
$k$--defective surface with reducible general tangential $k$--contact locus, hence 
$X$ sits in a cone of dimension $2k+1$, and not smaller, over a curve. \par

\end{itemize}\end{prop}

The previous results finish the classification in this case.

\begin{rem}\label {piudifett}\rm (1) The minimally $k$--defective threefolds
occurring in  Propositions \ref {proj}, \ref {proj'} and \ref {proj''} can be
$h$--defective for some $h>k$. \par

To be more precise, consider a 
threefold as in Proposition \ref {proj} and let $h>k$. We know that $X\subset \p^
r$ lies on a cone with vertex a subspace $\Pi$ of dimension $k-1$ over a surface
$S$. By applying Proposition \ref {conidif}, we see that $s^ {(h)}(X)= s^
{(h)}(Y)+k$. For example, if $S$ is not $h$--defective, then $s^
{(h)}(X)=\min\{r,3h+k+2\}$, so that $X$ is $h$--defective if and only if
$3h+k+2<r$. If $S$ is defective, $X$ is defective even for higher values of $h$.
We leave the details to the reader.

Similarly, if $X$ is as in Proposition \ref {proj'} [resp. Proposition 
\ref {proj''}] a similar argument shows that, for any $h>k$, one has $s^
{(h)}(X)=\min\{r,2h+2k+2\}$ [resp. $s^ {(h)}(X)=\min\{r,2h+2k+1\}$], so that $X$
is $h$--defective if and only if $2h+2k+2<r$ [resp. $2h+2k+1<r$]. \medskip

\noindent (2) A threefold as in Proposition \ref {proj} can be smooth. This can be
proved with a construction analogous to the one proposed in Remark \ref
{liscezza}, part (2). Similarly one can prove the existence of smooth threefolds
as in Proposition \ref {proj'}. Take indeed a smooth scroll surface $Y$ over a
curve $C$ spanning a $\p^{2k}$. Embed $C$ in a $\p^{r-2k-1}$. Then join every
point of $C$ with the corrresponding line of $Y$. The resulting scroll threefold
$X$ is smooth. An analogous construction works for producing smooth threefolds as
in Proposition \ref {proj''}.

\end{rem}

\section{THE IRREDUCIBLE CURVILINEAR CASE} \label {irredcurv}

Let $X\subset \pp^r$ be an irreducible, non--degenerate, projective,
minimally $k$--defective threefold with $k\geq 2$. In this section we
examine the case
$\gamma_k(X)=1$ and the general element of the family $\mathcal G$ is
irreducible.\par

By Proposition \ref {reduccontact} we know
that the general tangential
projection $X_{k-1}$ is a $1$--defective threefold whose general contact
locus is an
irreducible curve. From Remarkm \ref {remscorza} it follows that
$X_{k-1}$ as in cases (4) or (5) of Theorem \ref {scorza}.
In any event, we have that $X_{k-1}$ sits in $\pp^s$, $s=7,8,9$, hence:

$$ r\in
\{4k+3,4k+4,4k+5\}.$$

Now we can start our classification. The first step it to show that $X$ is as described in Proposition \ref {rnc}.

\begin{prop} Let $X\subset \pp^r$ be an irreducible, non--degenerate,
projective,
minimally $k$--defective threefold with $k\geq 1$. Assume
$\gamma_k=1$ and the general element $\Gamma$ of the family $\mathcal
G$ irreducible.
Then $\Gamma$ is a
rational normal curve in $\pp^{2k}$. Hence  $X$ is rationally
connected and therefore {\rm regular}, i.e. any desingularization $X'$ of $X$ has
$h^1(X',{\mathcal O}_{X'})=0$.\end{prop}

 \begin{proof} The assertion holds for $k=1$ (see Remark \ref {remscorza}). So we
may assume $k>1$ and proceed by induction on $k$.\par

Let $p\in X$ be a general point, let $\tau:=\tau_1:X\map X_1$ be the
tangential projection from
$T_{X,p}$ and let $\Gamma_1$ be the general tangential $(k-1)$--contact
locus of $X_1$. By
induction, $\Gamma_1$ is a rational normal curve in $\pp^{2k-2}$ and
$\Gamma$ maps onto
$\Gamma_1$ via $\tau$ (see the proof of Proposition \ref{reduccontact}).\par

One has $\dim (<\Gamma>)=2k$. Indeed by \eqref {terracini'} we have
$h(\Gamma)\leq 2k+1$, moreover the center of the projection $\tau$,
that is the  tangent space
$T_{X,p}$, meets $<\Gamma>$ at
least in the tangent line to $\Gamma$ at $p$ and the image of $\Gamma$
via $\tau$ spans a
$\pp^{2k-2}$. It follows also that, in particular, $<\Gamma>$ meets
$T_{X,p}$ exactly
along the aforementioned tangent line, so that $\tau_{|\Gamma}$ is a
general tangential
projection of $\Gamma$.\par

Now we claim that:\pas
\begin{itemize}

\item[(i)] $\Gamma$ is not tangentially degenerate;\pas

\item [(ii)] if $x\in \Gamma$ is a general point, the projection $\Gamma'$ of
$\Gamma$ from $x$ is also
not tangentially degenerate.
\end{itemize}

The assertions in the statement about $\Gamma$ follow from (i) and (ii). Indeed, by (ii) and
Lemma \ref {tanproj}, $\tau_{|\Gamma}$ is birational to its image
$\Gamma_1$, hence $\Gamma$ is
rational. Furthemore, since $\deg(\Gamma_1)=2k-2$ by induction, $<\Gamma>$ meets $T_{X,p}$ exactly
along $T_{C,p}$, as we saw, and $T_{C,p}\cap C=\{p\}$ by (i), we have $\deg(\Gamma)=2k$.

To prove (i) and (ii), we notice that $\mathcal G$ is a family of curves
of dimension $2(k+1)$. This is an easy count of parameters wich follows from the fact that there is a unique curve of $\mathcal G$ containing $k+1$ general points of $X$. Let
$p_1,\dots ,p_k$ be general points of $X$ and let $\mathcal G'$ be the
$2$--dimensional family of curves of $\mathcal G$ passing through
$p_1,\dots ,p_k$. We claim that:\pas

\begin{itemize}

\item [(iii)] the tangent lines to the curves of $\mathcal G'$ at $P_i$ fill up
the whole tangent space
$T_{X,p_i}$, for all $i=1,\dots ,k$.\par
\end{itemize}

Indeed, this is true for $k=1$ (see Remark \ref {remscorza}). Then proceed by
induction. The curves in
$\mathcal G'$ are mapped via $\tau_1=\tau_{X,p_1}$ to the contact
curves on $X_1$ through
$\tau_1(p_i)$, $i=2,\dots ,k$. Since the differential of $\tau_1$ is an
isomorphism at $p_2,\dots , p_k$,
which are general points, by induction the tangent lines to the curves
of $\mathcal G'$ at $p_i$ fill up the tangent space $T_{X,p_i}$,
for all $i=2,\dots ,k$. Arguing by monodromy (see \S {monodromy} and the proof of Proposition \ref{mon}) the same is true for $i=1$.\par

Now we claim that (iii) implies (i). Indeed, if $\Gamma$ were
tangentially degenerate, then, by (iii),
$X$ itself would be $2$-tangentially degenerate. By Proposition \ref
{tandeg}, $X$, which does not lie
in $\pp^4$, would be a scroll on a curve. Then also $X_{k-1}$ would be a
scroll, which is not the case (see Remark \ref{remscorza}).\par

Furthermore (iii) also implies (ii). In fact if we apply the same
argument as above to the projection
$X'$ of $X$ from a general point $x\in X$, we conclude that $X'$ would
be a scroll. Let $\Pi$ be a general plane of  $X'$. Then either $\Pi$ pulls back to a plane on $X$ or it pulls back to a quadric $Q$ containg $x$. In the former case $X$ would be a scroll, which, as we saw, is impossible. In the latter case $X$ is swept out by a family
$\mathcal Q$ of quadrics such that
through two general points of $X$ there is a quadric in $\mathcal Q$
containing them. In this case $X$ would be
$1$-defective (see Proposition \ref {rnc}), which is against
the minimality assumption.\end{proof}

To go on with the classification, we need the following:

\begin{ex} \label  {ex'}\rm  Let $Y\subset \p^ {k+2}$ be a threefold of minimal degree $k$ and let $X\subset \pp^{4k+5}$   be its double embedding. 
In \cite{WDV}, Example 4.5, it is proved that $X$ is
$(k+1)$--defective and $s^{(k+1)}(X)\leq 4k+4$. One can prove that  $X$ is actually
minimally $k$--defective,  that
$s^{(k)}(X)=4k+2$ and $s^{(k+1)}(X)= 4k+4$, so that $\delta_k(X)=1$, $n_k=1$.
Actually we will see that
$\gamma_k(X)=1$ and the general tangential $k$--contact locus is an irreducible
curve.\par

First we observe that $X$ is not $h$--weakly defective, hence not $h$--defective,
for any $h\leq k-1$. Indeed, if $p_0,...,p_h\in X$ are general points, there are
hyperplane sections of $X$ tangent at  $p_0,...,p_h$ and having isolated
singularities at  $p_0,...,p_h$. Take, for instance the hyperplane  sections
corresponding to sections of $Y$ with a  general quadric  cone with vertex along
the span of  the points corresponding to $p_0,...,p_h$. By the way, the same
argument applied for for $h=k$, shows  that there are
hyyerplane sections of $X$ tangent at  $p_0,...,p_k$ and having an irreducible
curve of singular  points containing  $p_0,...,p_k$, i.e. the rational normal
curve of degree $2k$ which  is the immage of the intersection of $Y$ with the
span of the points corresponding to $p_0,...,p_k$.\par

Now suppose $X$ is not $k$--defective. Then, if $p_0,...,p_k\in X$ are general
points, we would have $\dim(T_{X,p_0,...,p_k)})=4k+3$. However, if $p\in X$ is
another general point, one has $\dim(T_{X,p_0,...,p_k,p})\leq 4k+4$. This
would imply that the general tangent space $T_{X,p}$ to $X$ intersects
$T_{X,p_0,...,p_k}$ in dimension at least $2$. But this would force $X$ to span a
$\p^{4k+4}$ (see Proposition \ref {varincon}), a contradiction. Hence $X$
is minimally $k$--defective and therefore  $s^{(k)}(X)=4k+2$, $\delta_k(X)=1$,
$n_k=1$ (see Proposition \ref {difetti}). Then the same argument implies that 
$s^{(k+1)}(X)\geq 4k+4$ and therefore $s^{(k+1)}(X)= 4k+4$.\par

Of course also a projection of $X$ to $\p^{4k+3+i}$, $i=0,1$, is minimally
$k$--defective.

\end{ex}

To apply an inductive argument we will need the following:

\begin{prop}\label{4k+5} Let $X\subset \pp^r$ be an irreducible,
non--degenerate,
projective, regular threefold. Let $k\geq 2$ and assume that a general
tangential projection $X_1$ of $X$ is the $2$-uple embedding of a
minimal threefold in
$\pp^{k+1}$. Then $X$ is the $2$-uple embedding of a minimal threefold in
$\pp^{k+2}$.\end{prop}

\begin{proof} Let $\tau:=\tau_1: X\map X_1$ be the general tangential
projection from $T_{X,p}$.\par

Let $Z$ [resp. $Z_1$] be a desingularization of $X$ [resp. $X_1$] and
let $\mathcal H$
[resp. $\mathcal H_1$] be the pull--back on $Z$ [resp. $Z_1$] of the
hyperplane linear system of $X$
[resp. $X_1$]. The tangential
projection $\tau$ induces a rational map $\phi:Z \map Z_1$. Notice that:

\begin{equation} \label{pullback} \phi^*(\mathcal H_1)=\mathcal
H(-2p).\end{equation}

By the hypothesis $X_1$ spans a $\pp^{4k+1}$, so $r=4k+5$. Furthermore $X_1$ is
linearly normal, namely $\mathcal H_1$ is a complete linear system.
By (\ref {pullback}) also  $\mathcal H$ is complete, i.e.
$X$ is also linearly normal.\par

By assumption, $\mathcal H_1$ is the double of a linear system $\mathcal
L_1$ and the
associated map $\phi_{\mathcal L_1}: Z_1\map Y_1\subset\pp^{k+1}$ is such
that $Y_1$ is a
threefold of minimal degree. The linear system $\mathcal L_1$ pulls back
via $\phi$
to a linear system $\mathcal L_p$ on $Z$. Of course $\phi_{\mathcal
L_p}=\phi_{\mathcal L_1}\circ \phi:
Z\map  Y_1\subset\pp^{k+1}$ and therefore $\dim({\mathcal L_p})=k+1$. The
system $\mathcal
L_p$ depends on $p$, but, since $X$ is regular, as $p$ moves on $X$ all
the systems $\mathcal L_p$'s are subsystems of the same complete linear
system $\mathcal L$ on $Z$.\par

We notice that:

\begin{equation}\label{prima} 2\mathcal L_p=\phi^*(\mathcal
H_1)=\mathcal H(-2p).\end{equation}
Hence
all divisors in $\mathcal L_p$ contain $p$ and therefore we have

\begin{equation}\label{seconda}\mathcal L_p\subseteq\mathcal
L(-p)\end{equation}

As an imediate consequence of (\ref{prima}) and (\ref {seconda}) we have
the linear equivalence
relation $2\mathcal L\equiv \mathcal H$. Since $\mathcal H$ is complete,
we have $2\mathcal
L\subseteq \mathcal H$.\par

Moreover, by (\ref {seconda}), one has
$m:=\dim(\mathcal L)\geq k+2$. We call $Y\subset\pp^m$ the image of $Z$
via the map $\phi_{\mathcal
L}$ associated to $\mathcal L$.\par

If $m>k+2$, then $\dim(\mathcal H)\geq\dim(2\mathcal
L)>4k+5$ (see Theorem \ref {Hilb}), a contradiction. Similarly we
find a contradiction if
$\dim(\mathcal H)>\dim(2\mathcal L)$. Hence $\mathcal H=2\mathcal L$ and
$m=k+2$. Moreover for $p$ general, $\mathcal L_p=\mathcal L(-p)$. It
follows that $X$ is the
$2$-uple embedding of $Y$ and $Y\subset \pp^{k+2}$ projects from a
general point on it to
a minimal threefold of $\pp^{k+1}$. Hence $Y$ is also of minimal
degree.\end{proof}

\begin{rem}\label{4k+4-1} In the previous setting, assume that a general
tangential
projection $X_1$ of $X$ is a projection of the $2$-uple embedding $Y'$
of a threefold $Y$ of minimal
degree in $\pp^{k+1}$, from a space not intersecting $Y'$. Then $X$ is
the projection of the
$2$-uple embedding of a threefold of minimal degree in $\pp^{k+2}$, from
a space not intersecting
it.\medskip\rm

Indeed, by replacing, if necessary, $X_1$ by its linearly normal
embedding, the previous
argument shows that the linearly normal embedding of $X$ is the
$2$--uple embedding of a minimal threefold in $\pp^{k+2}$.  \end{rem}

A completely similar argument holds when $X_1$ is a projection of the
$2$-uple
embedding of a minimal threefold from a point on it. We skip the details.

\begin{prop}\label{4k+4-2} Let $X\subset \pp^r$ be an irreducible,
non--degenerate,
projective, regular threefold. Let $k\geq 2$ and assume that a general
tangential projection $X_1$ of $X$ is the projection of the $2$-uple
embedding $Y'$ of a
threefold $Y$ of minimal degree in $\pp^{k+1}$ from a space which meets
$Y'$ at a point. Then $X$ is
the projection of the $2$--uple embedding of a threefold of minimal
degree in $\pp^{k+2}$ from a
space which meets it at a point.\end{prop}

The case of projections from more than one point on $X$ is contained
in the analysis below.
Again, in order to apply an inductive argument, we will need the following:

\begin{prop}\label{PxP} Let $X\subset \pp^r$ be an irreducible,
non--degenerate,
projective, regular threefold. Let $k\geq 2$ and assume that a general
tangential projection $X_1$ of $X$ is linearly normal and contained in
the Segre embedding of
$\pp^{k}\times\pp^{k}$. Suppose that each of the two projections of $X_1$
to $\pp^k$ spans
$\pp^k$. Then $X$ is linearly normal and contained in the Segre embedding of
$\pp^{k+1}\times\pp^{k+1}$. Moreover each of the two projections of $X$
spans
$\pp^{k+1}$.\end{prop}

\begin{proof} Let $\tau:=\tau_1: X\map X_1$ be the general tangential
projection from $T_{X,P}$.\par

Let $Z$ [resp. $Z_1$] be a desingularization of $X$ [resp. $X_1$] and
let $\mathcal H$
[resp. $\mathcal H_1$] be the pull--back on $Z$ [resp. $Z_1$] of the
hyperplane linear system of $X$
[resp. $X_1$]. The tangential
projection $\tau$ induces a rational map $\phi:Z \map Z_1$. Equation
(\ref {pullback}) still
holds.  Then, as in the proof of Proposition \ref {4k+5}, one sees that
$X$ is linearly normal, i.e.
$\mathcal H$ is complete. \par

We have two linear systems $\mathcal A_1,\mathcal B_1$ on $Z_1$ which
come by pulling back the
linear system $|\mathcal O_{\pp^k}(1)|$ via the projections of
$\pp^{k}\times\pp^{k}$ to the
two factors. One has $\mathcal H_1=\mathcal A_1+\mathcal B_1$ (see \S \ref
{linsyst}.\par

The linear systems $\mathcal A_1,\mathcal B_1$
pull--back via $\phi$ to linear systems $\mathcal A_p, \mathcal B_p$, on
$Z$. By the regularity assumption on $X$, as $p$ varies $\mathcal
A_p, \mathcal B_p$ vary inside two complete linear systems $\mathcal A$,
$\mathcal B$
on $Z$. From the relation:

\begin{equation} \label {alpha}\mathcal A_p+\mathcal B_p=\phi^*(\mathcal
H_1)=\mathcal
H(-2p)\end{equation}
we see that every divisor in $\mathcal A_p+\mathcal B_p$ vanishes doubly
at $p$, and therefore either:
\begin{equation}\label{beta}\mathcal A_p\subseteq \mathcal A(-p),
\mathcal B_p\subseteq \mathcal
B(-p)\end{equation}
or, say:
\begin{equation}\label{beta'}\mathcal B_p\subseteq \mathcal
B(-2p)\end{equation}
In any case we have the equivalence relation $\mathcal A+\mathcal B\equiv
\mathcal H$. Since
$\mathcal H$ is complete, then

\begin{equation}\label {delta}\mathcal A+\mathcal B\subseteq \mathcal
H.\end{equation}

In case (\ref{beta}) holds, we have:

\begin{equation} \label {gamma}\mathcal A_p+\mathcal B_p\subseteq
\mathcal A(-p)+ \mathcal
B(-p)\subseteq \mathcal H(-2p)\end{equation}
 Then, by (\ref {alpha}), equality has to hold everywhere in
(\ref{gamma}) and therefore in
(\ref{beta}) and (\ref {delta}). This yields $\dim(\mathcal
A)=\dim(\mathcal B)=k+1$, proving the
assertion.\par

If (\ref{beta'}) holds, one has:

\begin{equation} \label {gamma'}\mathcal A_p+\mathcal B_p\subseteq
\mathcal A+ \mathcal
B(-2p)\subseteq \mathcal H(-2p)\end{equation}

 Again, by (\ref {alpha}), equality has to hold everywhere in
(\ref{gamma'}). From (\ref {delta}) we also
have  $\mathcal A(-p)+ \mathcal
B(-p)\subseteq \mathcal H(-2p)$ which is clearly incompatible with
$\mathcal H(-2p)=\mathcal A+
\mathcal B(-2p)$.\end{proof}

We can now get the following partial classification result:

\begin{thm}\label{reducib} Let $X\subset \pp^r$ be an irreducible,
non--degenerate,
projective, minimally $k$--defective threefold with $k\geq 2$. Assume
$\gamma_k=1$ and irreducible general tangential $k$--contact locus.
Then  one of the following holds:\pas

\begin{itemize}

\item [(1)] $r=4k+5$ and $X$ is the $2$-uple
embedding of a threefold of minimal degree in $\pp^{k+2}$;\pas

\item [(2)] $r=4k+4$ and $X$ is the
projection of the $2$-uple embedding of a threefold of minimal degree in
$\pp^{k+2}$ from a
point $p\in \pp^{4k+5}$; \pas

\item [(3)] $r=4k+3$ and either $X$ is the projection of the
$2$-uple embedding of a threefold of minimal degree in $\pp^{k+2}$ from
a line $\ell\subset
\pp^{4k+5}$, or $X$ is linearly normal, it is contained in the intersection of a
space of dimension $4k+3$ with the Segre embedding of  $\pp^{k+1}\times\pp^{k+1}$
in
$\pp^{k^2+4k+3}$ and each of the two projections of $X$ on the two factors 
of  $\pp^{k+1}\times\pp^{k+1}$ spans $\pp^{k+1}$. In this former case each of the
linear systems $\mathcal A_i$, $i=1,2$, on $X$
corresponding to the two projections $\phi_i: X\to \pp^{k+1}$, is base
point free and the
general surface $A_i\in \mathcal A_i$ is irreducible.\end{itemize} \end{thm}

\begin{proof} We know that the possibilities for $r$ are the ones listed in
statement.\par

Assume that $r=4k+5$. Then, by Theorem \ref {scorza}, $X_{k-1}$ is the
$2$-uple embedding of $\pp^3$. The conclusion follows by induction from
Proposition
\ref{4k+5}.\par

If $r=4k+4$, by Theorem \ref {scorza}, $X_{k-1}$ is a projection of the
$2$-uple embedding of $\pp^3$. The conclusion follows by induction from
Remark \ref{4k+4-1}
and Proposition \ref{4k+4-2}.\par

If $r=4k+3$, then by Proposition \ref {reduccontact}, by Theorem \ref {scorza} and
Remark \ref {remscorza}, either $X_{k-1}$ is some projection of the $2$-uple
embedding of
$\pp^3$ or $X_{k-1}$ is a
hyperplane section of $\pp^2\times\pp^2$. In the former case the
conclusion follows again by induction from Remark \ref{4k+4-1}
and Proposition \ref{4k+4-2}. In the
latter case, $X_{k-1}$ is linearly normal and the two projections of
$X_{k-1}$ to $\pp^2$ are both surjective. The conclusion follows by induction from
Proposition \ref{PxP}.
\end{proof}

The last part of case (3) of Theorem \ref {reducib} deserves more
attention. The analysis is based on the following proposition:

\begin{prop}\label{product} Let $Z$ be an irreducible, smooth, regular
threefold and let
$\mathcal A_i$, $i=1,2$, be two distinct base point free linear systems
of dimension $k+1\geq
3$ on $Z$ such that a general surface $A_i\in \mathcal A_i$ is smooth
and irreducible and
the minimal sum $\mathcal A_1+\mathcal A_2$ is birational on $Z$. Then
the minimal sum $\mathcal
A_1+\mathcal A_2$ has dimension at least $4k+3$.

 Furthermore if the minimum is attained then a general curve
$C$ in the class $A_1\cdot A_2$ is smooth, irreducible and rational and
the linear system
$\mathcal A_{i|C}$ is the complete $g^k_k$ on $C$. Moreover, either:\pas

\begin{itemize}

\item [(a)] $\mathcal A_i$, $i=1,2$, is
complete, or\pas
\item [(b)] $A_1\equiv A_2$, and $\dim(|A_i|)=k+2$, $i=1,2$, or\pas
\item [(c)] $\mathcal A_1$ is complete, whereas $\dim(|A_2|)=k+2$ and there is
an effective, non--zero divisor
$E$ on $X$, such that $h^0(Z,\mathcal O_Z(E))=1$, $A_2-A_1\equiv E$ and
$\mathcal O_{A_2}(E)=\mathcal
O_{A_2}$.\par
\end{itemize}

If, in addition, the linear systems $\mathcal A_i$,
$i=1,2$, are very big, then also the general curves $C_i$ in the classes
$A_i^2$ are smooth,
irreducible and rational. Moreover:\pas

\begin{itemize}

\item [(a')] in case (a) above  the linear system $\mathcal
A_{i|C_i}$ is the complete $g^{k-1}_{k-1}$ on $C_i$, $i=1,2$;\pas

\item [(b')] in case (b) above  the linear system $\mathcal
A_{i|C_i}$ is the complete $g^{k}_{k}$ on $C_i$, $i=1,2$;\pas

\item [(c')] in case (c) above the linear system $\mathcal
A_{1|C_1}$ is the complete $g^{k-1}_{k-1}$ on $C_1$ whereaes the linear
system $\mathcal
A_{2|C_2}$ is the complete $g^{k}_{k}$ on $C_2$.\end{itemize}\end{prop}

\begin{proof} The curve $C$ is smooth by Bertini's theorem. We claim
that $C$ is also irreducible.
Otherwise, by Bertini's theorem, the linear system $\mathcal A_{1|A_2}$ would be
composite with a pencil $\mathcal P$ and, if
$P\in \mathcal P$ is a general curve, we would have $P\cdot A_1=0$.
Similarly we would have $P\cdot
A_2=0$ and therefore each curve $P$ would be contracted by the map  $\phi_{\mathcal
A_1+\mathcal A_2}$, against the
hypothesis that $\mathcal A_1+\mathcal A_2$ is birational.\par

Consider the linear systems $\mathcal A_{i|C}$, $i=1,2$. Since $\mathcal
A_1\not=\mathcal A_2$, one has:

\begin{equation}\label {kappa'} \dim(\mathcal A_{i|A_{j}})=k+1, i,j=1,2,
i\neq j\end{equation}
 and therefore:

\begin{equation} \label{kappa}
\dim(\mathcal A_{i|C})=k, i=1,2.\end{equation}
By Lemma \ref {hopf} we have $\dim(\mathcal
A_{1|C}+\mathcal A_{2|C})\geq 2k$.\par

  Let $\dim(\mathcal A_{1|A_1}+\mathcal A_{2|A_1})=s$. By
looking at the exact sequence:

$$0\to H^0(Z,\mathcal O_Z(A_2))\to H^0(Z,\mathcal O_Z(A_1+A_2))\to
H^0(A_1,\mathcal
O_{A_1}(A_1+A_2))$$
we see that the image of the rightmost map contains the vector space
corresponding to $\mathcal
A_{1|A_1}+\mathcal A_{2|A_1}$. This implies that:

$$\dim(\mathcal A_1+\mathcal A_2)\geq s+1+\dim(\mathcal A_2)=s+k+2.$$
Similarly, by looking at the sequence:

$$0\to H^0(A_1,\mathcal O_{A_1}(A_1))\to H^0(A_1,\mathcal
O_{A_1}(A_1+A_2))\to H^0(C,\mathcal
O_{C}(A_1+A_2))$$
one obtains:

$$s+1\geq \dim(\mathcal A_{1|C}+\mathcal A_{2|C})+\dim(\mathcal
A_{1|A_1})+2\geq 3k+2$$
whence $\dim(\mathcal A_{1}+\mathcal A_{2})\geq 4k+3$ follows.\par

If the equality holds, then Lemma \ref {hopf} implies that $C$
is rational. Since, by the hypotheses and by Lemma \ref {big}, the
linear series $\mathcal
A_{1|C}+\mathcal A_{2|C}$ is birational, then Lemma \ref {hopf} again
implies that $\mathcal
A_{1|C}=\mathcal A_{2|C}$ and that this is a complete linear series.
Since $C$ is rational, by (\ref
{kappa}) we have $A_i\cdot C=k$.\par

Since $\mathcal
A_{1|C}=\mathcal A_{2|C}$ are complete, it follows that the systems
$\mathcal
A_{1|A_2}$ and $\mathcal A_{2|A_1}$ are also complete. If $h^0(Z,\mathcal
O_Z(A_1-A_2))=h^0(Z,\mathcal O_Z(A_2-A_1))=0$, then, by looking at the
exact sequence:

\begin{equation} \label {restrict} 0\to H^0(Z,\mathcal O_Z(A_i-A_j))\to
H^0(Z,\mathcal O_Z(A_i))\to
H^0(A_j,\mathcal O_{A_j}(A_i)), i,j=1,2, i\neq j \end{equation}
we see we are in case (a).\par

Suppose we are not in case (a) and assume that $h^0(Z,\mathcal
O_Z(A_2-A_1))>0$, so that there is an effective divisor $E$ which is
linearly equivalent to
$A_2-A_1$. Suppose $E$ is zero. Then $A_1\equiv A_2$. Moreover since
$\mathcal A_{2|A_1}$ is
complete, from the sequence (\ref {restrict}) with $i=2, j=1$,
we see that $h^0(Z,\mathcal O_Z(A_2))=\dim(\mathcal A_{2|A_1})+2=k+3$,
namely we are in case
(b).\par

If $E$ is non zero, then $h^0(Z,\mathcal
O_Z(A_1-A_2))=0$. Since $\mathcal
A_{1|A_2}$ is complete, the completeness of $\mathcal A_1$ follows from
the exact sequence (\ref
{restrict}) with $i=1, j=2$.
Moreover $E\cdot C=0$. Since the linear system $|C|$ on $A_2$ is base point free,
with positive self--intersection, we see that $h^0(A_2,\mathcal O_{A_2}(E))=1$.
From the exact sequence:

$$0\to H^0(Z, \mathcal O_Z(-A_1))\to H^0(Z,\mathcal O_Z(E))\to
H^0(A_2,O_{A_2}(E))$$ we deduce that $h^0(Z,\mathcal \mathcal O_{Z}(E))=1$. Now
look at the sequence (\ref {restrict}) with $i=2,
j=1$.
Since, as we saw, $h ^0(A_1,\mathcal
O_{A_1}(A_2))=\dim(\mathcal A_{2|A_1})+1=k+2$, we have $ h
^0(Z,\mathcal O_{Z}(A_2))=k+3$, i.e.
$\dim(|A_2|)=k+2$. Now, inside $|A_2|$ we have the two linear systems
$E+\mathcal A_1$ and $\mathcal
A_2$. They intersect along a $k$--dimensional linear system. Hence every
section of $H ^0(Z,\mathcal
O_{Z}(A_2))$ determining a divisor of $\mathcal A_2$ is of the form
$fe+g$, with $f\in H ^0(Z,\mathcal
O_{Z}(A_1))$ variable, $e\in H ^0(Z,\mathcal
O_{Z}(E))$ and $g\in H ^0(Z,\mathcal
O_{Z}(A_2))$ fixed defining a divisor $B\mathcal A_2$ not containing $E$. Hence
every solution of the system
$e=g=0$, i.e. any point in $E\cap B$, should
be a base point of $\mathcal A_2$. Since $\mathcal A_2$ is
base point free, we have that $E\cap B=\emptyset$. Thus we are in case
(c).\par

Suppose now $\mathcal A_i$,
is very big. Arguing as at the beginning of the proof, we see that the
general curve
$C_i$ in the class $A_i^2$ is smooth and irreducible. By what we proved
already, we have the exact
sequence:

$$0\to \mathcal O_{A_1}(-A_1)\to \mathcal O_{A_1}(A_2-A_1)\to \mathcal
O_{C}(A_2-A_1)\simeq
\mathcal O_{C}\to 0$$
from which we deduce $h^0(A_1,O_{A_1}(A_2-A_1))=1,
h^1(A_1,O_{A_1}(A_2-A_1))=0$. Indeed, by Lemma
\ref {big}, $\mathcal O_{A_1}(A_1)$ is big and nef and therefore
Kawamata-Viehweg  vanishing theorem
says that $h^i(A_1,\mathcal O_{A_1}(-A_1))=0$, $0\leq i\leq 2$.\par

Now look at the sequence:

$$0\to \mathcal O_{A_1}(A_2-A_1)\to \mathcal O_{A_1}(C)\to \mathcal
O_{C_1}(C)\to 0.$$
Since $\mathcal A_{2|A_1}\subset |\mathcal O_{A_1}(C)|$, by (\ref
{kappa'}) we have
$h^0(A_1,O_{A_1}(C))\geq k+2$. It follows that $h^0(C_1,O_{C_1}(C))\geq
k+1$. Since
$\deg(O_{C_1}(C))=k$, we deduce that $C_1$ is rational.\par

Suppose we are in case (a). Since
$Z$ is regular and $\mathcal A_1$ is complete, we have:

$$h^0(A_1,\mathcal
O_{A_1}(C_1))=h^0(A_1,\mathcal O_{A_1}(A_1))=\dim(\mathcal A_1)=k+1.$$
Finally look at the sequence:

$$0\to \mathcal O_{A_1}\to \mathcal O_{A_1}(C_1)\to \mathcal
O_{C_1}(C_1)\to 0.$$
Clearly $A_1$ is a rational surface, so we have $h^1(A_1,\mathcal
O_{A_1})=0$.
Thus one  has $h^0(C_1,\mathcal O_{C_1}(C_1))=k$. This proves the
assertion for $\mathcal A_{1|C_1}$.
The same for $C_2$. Thus we are in case (a').\par

The analysis in case (b) and (c) is the same, leading to (b') and (c')
respectively. \end{proof}

\begin{cor} \label {finalcor} Let $X$ be an irreducible, non--degenerate
threefold in the Segre embedding of
$\pp^{k+1}\times\pp^{k+1}$, $k\geq 2$, which does not lie in the
$2$--uple embedding of
$\pp^{k+1}$. Assume that each of the two projections of $X$ to
$\pp^{k+1}$ spans $\pp^{k+1}$.
Then $X$ spans a space of dimension at least $4k+3$.\par

Furthermore, if $X$ spans a  $\pp^{4k+3}$, then given $k+1$ general
points of $X$, there is
a rational normal curve $C$ of degree $2k$ on $X$ containing the given
points, and $X$ is minimally
$k$--defective and  $s^{(k)}(X)=4k+2$, $\delta_k(X)=1$, $n_k(X)=2$.

Finally, if the two projections of $X$ to $\pp^{k+1}$ are generically
finite to
their images, then either:\pas

\begin{itemize}

\item [(a)] they both map birationally $X$ to rational normal scrolls, and then
the degree of $X$ is $8k-2$,
or\pas

\item [(b)] they both map birationally $X$ to projections of rational normal
scrolls in $\pp^{k+2}$ from a
point, and then the degree of $X$ is $8k$, or\pas

\item [(c)] one of them maps birationally $X$ to a rational normal scroll and
the other maps birationally $X$
to the projection of a rational normal scroll in $\pp^{k+2}$ from a
point, and then the degree of $X$
is $8k-1$.\end{itemize}\end{cor}

\begin{proof} The first assertion follows from Proposition \ref{product}
applied to a
desingularization $Z$ of $X$. Consider the two linear systems $\mathcal
A_i$, $i=1,2$, on $Z$
corresponding to the two projections $X\to \pp^{k+1}$. 
 Notice that the general surface of
$\mathcal A_i$, $i=1,2$, is irreducible (see Proposition \ref {reducib}) and smooth
by Bertini's theorem, since $\mathcal A_i$ is base point free.\par

Suppose $X$ spans a  $\pp^{4k+3}$. Given $k+1$ general points
$p_0,\dots,p_k$ on $Z$, let $A_i\in \mathcal A_i$ be the surface
containing $p_0,\dots,p_k$.
Then, according to Proposition \ref{product}, the image in $X$ of the
curve $C=A_1\cdot A_2$ is a
rational normal curve of degree $2k$. The $k$--defectivity of $X$ follows from
Proposition
\ref {rnc}. \par

Let us prove that $X$ is minimally $k$--defective. We first claim that $X$ is not
$1$--defective. Assume, by contradiction, that $X$ is $1$-defective.  Then, by Theorem
\ref {scorza}, the only possibilities are that $X$ is either a cone, or $X$ sits
in a $4$--dimensional cone over a curve. Notice that
the Segre embedding of
$\pp^{k+1}\times\pp^{k+1}$ is swept out by two $(k+1)$--dimensional families of
$\p^{k+1}$'s and that each line on the Segre embedding
of $\pp^{k+1}\times\pp^{k+1}$ is contained in a
$\p^{k+1}$ of either one of these two families. As a consequence, each
irreducible cone contained in the Segre embedding of
$\pp^{k+1}\times\pp^{k+1}$ is contained in a $\pp^{k+1}$. This implies that $X$
cannot be a cone, since it spans a $\p^{4k+3}$. Suppose $X$ sits
in a $4$--dimensional cone over a curve. Let $p\in X$ be a general point and
consider the general tangential projection $\tau$ of $\pp^{k+1}\times\pp^{k+1}$ 
from $p$. As we saw in Example \ref {segre}, $\pp^{k+1}\times\pp^{k+1}$  projects
onto
$\pp^{k}\times\pp^{k}$. Let $X'$ be the image of $X$ via $\tau$. Notice that $X'$
is a projection of the image $X_1$ of the tangential projection of $X$ from
$T_{X,p}$. Thus $X'$ is a cone over a surface. By the above argument, $X'$
would span at most a $\p^k$. Thus $4k+3=\dim(<X>)=\dim(<X'>)+2k+3\leq 3k+3$, a contradiction.\par

Let us consider again a general point $p\in X$ and the general projection
$\tau_{p}$ from $T_{X,p}$. Since $X$ is not $1$--defective, its image is a
threefold spanning a
$\p^{4k-1}$. Consider again a desingularization $Z$ of $X$. We abuse notation
and denote by
$p$ the point of $Z$ corresponding to $p\in X$. Let $\mathcal H$ be the linear
system on $Z$ corresponding to the hyperplane section system on $X$.
We have $\mathcal A_1(-p)+\mathcal
A_2(-p)\subseteq \mathcal H(-2p)$. On the other hand, since $\dim(\mathcal
A_i(-p))=k$, $i=1,2$, by Proposition \ref
{product}, we have $4k-1=\dim(\mathcal H(-2p))\geq \dim(\mathcal A_1(-p)+\mathcal
A_2(-p))\geq 4k-1$. This proves that $\mathcal A_1(-p)+\mathcal
A_2(-p)=\mathcal H(-2p)$. On the other hand $\phi_{\mathcal H(-2p)}$ just maps $Z$
to $X_1$. By what we saw in \S \ref {linsyst}, we have that $X_1$ is a threefold
which sits in $\pp^{k}\times\pp^{k}$, spans a $\p^{4k-1}$ and 
each of the two projections of $X_1$ to
$\pp^{k}$ spans $\pp^{k}$. Thus, by arguing as above, we see that $X_1$ is not
$1$--defective. By iterating this argument, one proves that $X_i$ is not
$1$--defective for any $i=1,...,k-2$, hence that $X$ is minimally $k$--defective.
\par

Now, notice that the above argument proves that $X_{k-1}$ is the hyperplane
section of the Segre embedding of $\p^2\times \p^2$, hence 
$s^{(1)}(X_{k-1})=6$,
$\delta_1(X_{k-1})=1$, $n_1(X_{k-1})=2$. Thus
$s^{(k)}(X)=4k+2$, $\delta_k(X)=1$, $n_k(X)=2$.

Finally, the cases (a), (b) and (c) in the statement, correspond to the
homologous case in Proposition
\ref{product}. \end{proof}

The following proposition takes care of the cases left out in the
statement of Corollary \ref
{finalcor}.

\begin{prop}\label{veronese} Let $X$ be an irreducible, non--degenerate
threefold in the Segre embedding of
$\pp^{k+1}\times\pp^{k+1}$, $k\geq 2$, which lies in the $2$-uple
embedding of
$\pp^{k+1}$. Suppose the two projections of $X$ to $\pp^{k+1}$ span
$\pp^{k+1}$. Then these
projections coincide and are immersions of $X$ in $\pp^{k+1}$.
Furthermore the dimension $s$ of the
linear span of $X$ satisfies $s>4k+3$ unless:\pas

\begin{itemize}

\item [(a)] $s=4k+1$ and the image of $X$ in $\pp^{k+1}$ is a threefold $Y$ of
minimal degree;\pas

\item [(b)] $s=4k+2$ and the image of $X$ in $\pp^{k+1}$ is a threefold $Y$ of
degree $k$
with curve sections of arithmetic genus $1$;\pas

\item [(c)]
 $s=4k+3$ and the image of $X$ in $\pp^{k+1}$ is:\pas

\begin{itemize}

\item [$(c_1)$] either a threefold $Y$ of degree $k+1$ with curve sections of
arithmetic genus $2$ or\pas

\item [$(c_2)$] a
threefold $Y$ of degree $k$ which is the projection in $\pp^{k+1}$ of a
threefold of minimal degree in
$\pp^{k+2}$. In the latter case $Y$ is described in Lemma \ref {RNS}.\par
\end{itemize}
\end{itemize}

The threefolds in case (a) are $(k-1)$--defective, whereas the threefolds
in cases (b), (c) are
minimally $k$--defective. \end{prop}

\begin{proof} In the present situation it is clear that the two
projections of $X$ to $\pp^{k+1}$
 coincide and are immersions $\phi: X\to Y\subset \pp^{k+1}$ of $X$ in
$\pp^{k+1}$. Let $d$ be the degree of $Y$. Set $d=k-1+\iota$. Now we
make free and iterate use of
Theorem \ref {Hilb}. So we have $s\geq 4k+1+\iota$. If
$\iota=0$ we are in case (a). If
$\iota=1$ and $s=4k+2$  we are in case (b). If  $\iota=1$ and
$s=4k+3$, the curve sections of $Y$ cannot have arithmetic genus $1$,
otherwise $s=4k+2$. Then the
curve sections of $Y$ have arithmetic genus $0$, and we are in  case
$(c_2)$, as described in Lemma
\ref {RNS}. Finally if $\iota=2$ and $s=4k+3$ we are in case $(c_1)$. In
all other cases $s\geq
4k+4$.\par

The defectivity of the varieties in the list follows by Example \ref
{ex'}, (6), Example \ref
{exx} (3), (4), (6).\end{proof}

\begin{ex}\label {biproj} \rm  (1) It is clear that there are examples of 
smooth threefolds as in cases (1) and (2) of Theorem \ref {reducib}. It is also not
difficult to find examples of smooth threefolds as in case (3) of the samme
theorem, and specifically as in cases (a), (b), (c) of Corollary \ref
{finalcor}.\par

 For instance,  take a threefold $Y$ of minimal degree in $\pp^{k+2}$
and call $\mathcal H$
its hyperplane linear system. Put $\mathcal A_i=\mathcal H(-p_i)$,
$i=1,2$, with $p_1,p_2\in
\pp^{k+2}$ distinct points. Then $\mathcal A_1+\mathcal A_2$ embeds $Y$ in the
product
$\pp^{k+1}\times\pp^{k+1}$ and
also in $\pp^{4k+3}$, since it corresponds to a projection of the
$2$--uple embedding of $Y$
from the line $\ell=<p_1,p_2>$. Notice that these projections have degree
$d$ with $d=8k, 8k-1, 8k-2$,
according to the fact that $\ell$ does not intersect $Y$, meets $Y$ at
one point, or it is a secant of $Y$. \medskip

\noindent (2) There are examples of minimally $k$-defective threefolds, contained
in the intersection of the Segrre
embedding of $\pp^{k+1}\times\pp^{k+1}$ with a $\pp^{4k+3}$, for which one the
linear systems
$\mathcal A_i$ of Proposition \ref{product} is not very big.\par

This is the case of the threefold $X$ obtained by embedding
$\pp^1\times\pp^2$ in $\pp^{19}$ with the complete system $\mathcal
A$ of divisors of type $(1,3)$. Terracini pointed out in \cite{Terr3} that
$X$ is $4$--defective. Terracini's paper was reconsidered, from a modern
point of view, by Dionisi and Fontanari (\cite {DF}), who also proved that
the $X$ is not $3$--defective.
Let us show how this threefold fits in our classification.\par

The linear system $\mathcal A$ can be decomposed as the sum $\mathcal
A_1+\mathcal A_2$ with $\mathcal A_1$ of type $(1,1)$, $\mathcal
A_2$ of type $(0,2)$. Both $\mathcal A_1$, $\mathcal A_2$ send
$\pp^1\times\pp^2$ to $\pp^5$, so $X$ also sits in the product
$\pp^5\times\pp^5$, and the projections to the two factors are both
non--degenerate. Therefore $X$ fits in the last case (15) of our
classification. According to Proposition \ref{product}, a general curve
in the class $\mathcal A_1\cdot\mathcal A_2$ is rational: indeed it is
embedded as a rational normal quartic in the natural Segre embedding of
$\pp^1\times\pp^2$. \par

Observe, however, that $\mathcal A_2$ is 
not big: the second projection of $\pp^5\times\pp^5$ sends $X$ to a
surface. So we cannot apply to $X$ the last part of Proposition
\ref{product} as well as the last part of Corollary \ref{finalcor}. In
particular the degree of $X$ is $27<30=8k-2$ and $X$ is not the projection of
the $2$--uple embedding of a quartic threefold in $\pp^6$ from a line.\par

Other similar examples of this sort can be considered and are related to
the so--called phenomenon of {\it Grassmann defectivity} (see \cite {DF}).
We hope to come back on this subject in a future paper.
\end{ex}

\end{document}